\documentclass[12pt]{elsarticle}
\usepackage{amsmath,amssymb,epsfig,amsfonts,epsfig,graphicx}
\usepackage{stmaryrd} 
\usepackage{epsfig}
\usepackage{color}
\usepackage{bbold}

\newcommand{\cblue}{\color{black}}

\newcommand{\cn}{\color{black}}

\renewcommand{\phi}{\varphi}

\journal{Mathematics and Mechanics of Solids}

\newtheorem{Remark}{Remark}
\newtheorem{Lemma}{Lemma}
\newtheorem{Corollary}{Corollary}

\renewcommand{\rho}{\varrho}
\newcommand{\dx}{\,\mathrm{dx}}
\newcommand{\ddx}{\frac{\mathrm{d}}{\mathrm{dx}}}

\newcommand{\pat}{\partial_t}

\newcommand{\sa}{\sin(\alpha)}
\newcommand{\ca}{\cos(\alpha)}
\newcommand{\ta}{\tan(\alpha)}
\newcommand{\ssa}{\sin^2(\alpha)}
\newcommand{\cca}{\cos^2(\alpha)}
\newcommand{\shalf}{\sin^2\!\big(\frac{\alpha}{2}\big)}

\newcommand{\DDiv}{\mathrm{Div}}
\newcommand{\Wred}{W_\mathrm{red}}
\newcommand{\Ered}{E_\mathrm{red}}

\newcommand{\taured}{\tau_\mathrm{red}}
\newcommand{\uhom}{\overline{u}}
\newcommand{\cXD}{{\cal X}^D}
\newcommand{\XD}{X^D}
\newcommand{\cXCred}{{\cal X}_\mathrm{red}^C}
\newcommand{\XCred}{X_\mathrm{red}^C}

\newcommand{\nn}{\nonumber}
\newcommand{\io}{\int_0^1}

\newcommand{\cL}{{\cal{L}}}
\newcommand{\cM}{\cal{M}}
\newcommand{\cX}{{\cal{X}}}
\newcommand{\weakY}{\stackrel{Y}{\to}}

\newcommand{\R}{\mathbb{R}}
\newcommand{\N}{\mathbb N}
\newcommand{\Z}{\mathbb Z}
\newcommand{\tr}{\mathrm{tr}}

\newcommand{\SO}{\mathrm{SO}}

\newcommand{\DR}{\text{D}}

\newcommand{\Curl}{\text{Curl}}
\newcommand{\dev}{\text{dev}}
\newcommand{\sym}{\text{sym}}
\newcommand{\id}{\mathbb{1}}
\newcommand{\iprod}{\innerproduct}
\newcommand{\so}{\mathfrak{so}}
\newcommand{\innerproduct}[1]{\langle #1 \rangle}
\renewcommand{\skew}{\mathop{\mathrm{skew}}\nolimits}
\DeclareMathOperator*{\argmin}{argmin}
\DeclareMathOperator{\Sym}{sym}
\DeclareMathOperator{\skw}{skew}

\begin{document}

\begin{frontmatter}
\title{Simple shear in nonlinear Cosserat micropolar elasticity: Existence of
minimizers, numerical simulations and occurrence of microstructure}
\author{Thomas Blesgen}
\ead[1]{t.blesgen@th-bingen.de}
\address{Bingen University of Applied Sciences, Berlinstra{\ss}e 109,
D-55411 Bingen, Germany}
\author{Patrizio Neff}
\ead[2]{patrizio.neff@uni-due.de}
\address{University of Duisburg-Essen, Faculty of Mathematics,
Thea-Leymann-Stra{\ss}e 9, D-45127 Essen, Germany}
\cortext[2]{Corresponding author}
\cormark[2]

\begin{abstract}
Deformation microstructure is studied for a $1D$-shear problem in geometrically
nonlinear Cosserat elasticity. Microstructure solutions are described
analytically and numerically for zero characteristic length scale.
\end{abstract}

\end{frontmatter}
\vspace*{-3mm}
\begin{center}
\cblue Revised version. \cn
\end{center}


\section{Introduction}
\label{secintro}
This article studies the formation of microstructure due to simple shear
boundary conditions within a geometrically nonlinear Cosserat theory.

The Cosserat model is one of the best-known generalized continuum models
\cite{Capriz89}. It assumes that material points can undergo translation,
described by the standard deformation map $\varphi: \Omega\to\R^3$ and
independent rotations described by the orthogonal tensor field
$R: \Omega\to\SO(3)$, where $\Omega\subset\R^3$ describes the reference
configuration of the material. Therefore, the geometrically nonlinear Cosserat
model induces immediately the Lie-group structure on the configuration space
$\R^3\times\SO(3)$. 

\newpage
Both fields are coupled in the assumed elastic energy $W=W(\DR\varphi,R,\DR R)$
and the static Cosserat model appears as a two-field minimization problem which
is automatically geometrically nonlinear due to the presence of the non-abelian
rotation group $\SO(3)$. Material frame-indifference (objectivity) dictates
left-invariance of the Lagrangian $W$ under the action of $\SO(3)$ and material
symmetry (here isotropy) implies right-invariance under action of $\SO(3)$.

In the early 20th century the Cosserat brothers E. and F. Cosserat introduced
this model in its full geometrically nonlinear splendor
\cite{cosserat1909theorie} in a bold attempt to unify field theories embracing
mechanics, optics and electrodynamics through a common principal of least
action. They used the invariance of the energy under Euclidean transformations
\cite{cosserat1991note} to deduce the correct form of the energy
$W=W(R^T\DR\varphi, R^T\partial_xR,R^T\partial_yR,R^T\partial_zR)$ and to
derive the equations of balance of forces (variations w.r.t. the deformation
$\varphi$, the force-stress tensor may loose symmetry,
\cite{neff2008symmetric}), and balance of angular momentum
(variations w.r.t. rotations $R$).
The Cosserat brothers did not provide, however, any specific constitutive form
of the energy since they were not interested in specific applications.

\subsection{Three-dimensional geometrically nonlinear isotropic Cosserat model}
The underlying three-dimensional isotropic Cosserat model can be described in
terms of the standard deformation mapping $\varphi:\Omega\subset\R^3\to\R^3$
and an additional orthogonal microrotation tensor
$R:\Omega\subset\R^3\to \SO(3)$.

The goal here is to find a minimizer of the following isotropic
energy\footnote{The volumetric term $\frac{\lambda}{4}\big(
(\det\overline{U}-1)^2+(\frac{1}{\det\overline{U}}-1)^2\big)$ is independent of
the microrotation $R$ and polyconvex in $\DR\varphi$. It's quadratic
approximation is $\frac{\lambda}{2}\tr^2(\overline{U}-\id_3)$.}
\small
\begin{align}
E(\varphi,R) &= \int_\Omega\mu \big|\sym(\overline{U}-\id_3)\big|^2
+\mu_c\big|\text{skew}(\overline{U}-\id_3)\big|^2
+\frac{\lambda}{4}\big((\det\overline{U}-1)^2+(\frac{1}{\det\overline{U}}-1)^2
\big)\nn\\
&\quad\quad +\mu\frac{L_c^2}{2}\Big(a_1\big|\dev\, \sym\, R^T\Curl R\,
\big|^2+a_2\big|\text{skew}R^T\Curl R\big|^2+\frac{a_3}{3}\tr(R^T\Curl R)^2
\Big)\dx\nn\\
&= \int_{\Omega}W_{\text{mp}}(\overline{U})+W_{\text{disloc}}(R^T\Curl R)\dx
\to\min\quad\text{w.r.t.}\quad (\varphi,R)\,,\qquad
\overline{U}=R^T\DR\varphi\,.
\label{isotrop ener}
\end{align}
The problem will be supplemented by Dirichlet boundary conditions for the
deformation $\varphi$ and the microrotations $R$ can either be left free or
prescribed or connected to $\DR\varphi$ via the coupling condition
$\skew(\overline{U})|_{\Gamma}=0$ with $\Gamma=\partial\Omega$.
Here, $\mu>0$ is the standard elastic shear modulus, $\lambda$ the second
elastic Lam\'{e} parameter and $\mu_c\geq0$ is the so-called Cosserat couple
modulus; $a_1,a_2,a_3$ are non-dimensional non-negative weights and $L_c>0$ is
a characteristic length. The isotropic energy (\ref{isotrop ener}) is written in
terms of the nonsymmetric Biot type stretch tensor $\overline{U}=R^T\DR\varphi$
(first Cosserat deformation tensor, \cite{cosserat1909theorie}) and the
curvature measure $R^T\Curl R$. We call ${\bf{{\bf{\alpha}}}}:=R^T\Curl R$
the \textit{second order} dislocation density tensor,
\cite{birsan2016dislocation}.
Due to the orthogonality of $\dev\;\sym$, $\text{skew}$ and $\tr(.)\id$, the
curvature energy provides a complete control of 
\begin{align}
|{\bf{\alpha}}|^2=\big|R^T\Curl R\big|^2=|\Curl R|^2\qquad\text{provided}
\qquad a_1,a_2,a_3>0\,.
\end{align}
%
Using the result in \cite{neff2008curl},
\begin{align}
\big|\Curl R\,\big|^2_{\R^{3\times 3}}\geq c^+
|\DR R\,|^2_{\R^{3\times 3\times 3}}\,,
\end{align}
shows that the energy (\ref{isotrop ener}) controls $\DR R$ in
$L^2(\Omega,\R^{3\times 3\times 3})$. 

In this setting, the minimization problem is highly non-convex w.r.t.
$(\varphi,R)$. Existence of minimizers for (\ref{isotrop ener}) with $\mu_c>0$
has been shown first in \cite{neff2004existence}, see also
\cite{neff2015existence, neff2004existence, lankeit2017integrability,
birsan2016dislocation, mariano2009ground}, the partial
regularity for minimizers of a related problem is investigated in
\cite{li2022regularity, tel2019regularity}.
The Cosserat couple modulus $\mu_c$ controls the deviation of the microrotation
$R$ from the continuum rotation $\text{polar}(\DR \varphi)$ in the polar
decomposition of
$\DR\varphi=\text{polar}(\DR\varphi)\cdot\sqrt{\DR\varphi^T\DR \varphi}$,
cf. \cite{neff2014grioli}.

For $\mu_c\to \infty$ the constraint $R=\text{polar}(\DR\varphi)$ is generated
and the model would turn into a Toupin couple stress model.

Here we derive the three-dimensional Euler-Lagrange equations based on the
curvature expressed in the dislocation tensor ${\bf{\alpha}}=R^T\Curl R$.
We can write the bulk elastic energy as
\begin{align}
\label{Euler1}
E(\varphi,R)=\int_\Omega W_{\text{mp}}(\overline{U})
+W_{\text{disloc}}({\bf{\alpha}})\;{\rm{dx}}\,,\quad\quad
\overline{U}=R^T\DR \varphi\,,\quad\quad {\bf{\alpha}}=R^T\Curl R\,.
\end{align}
Taking variations of (\ref{Euler1}) w.r.t. the deformation $\varphi$ leads to
\begin{align}
\delta E{(\varphi,R)}\cdot\delta\varphi=\int_{\Omega} &
\iprod{\DR W_{\text{mp}}(\overline{U}),R^T\DR \delta\varphi}_{\R^{3\times 3}}
\;{\rm{dx}}=0\,,\qquad\forall\;\delta\varphi\in C_0^\infty(\Omega,\R^3)\\
\nn\Longleftrightarrow & \int_{\Omega}\iprod{R\,\DR
W_{\text{mp}}(\overline{U}),\DR\delta\varphi}_{\R^{3\times 3}}\;{\rm{dx}}
=\int_{\Omega}\iprod{\DDiv[R\cdot\DR W_{\text{mp}}(\overline{U})],
\delta\varphi}_{\R^3}\,{\rm{dx}}=0\,.
\end{align}
Taking variation w.r.t. $R\in\SO(3)$ results in (abbreviate $F:=\DR \varphi$)
\begin{align}
\label{Euler2}
\delta E{(\varphi,R)}\!\cdot\!\delta R &= \int_{\Omega}
\iprod{\DR W_{\text{mp}}(\overline{U}),\delta R^TF}
+\iprod{\DR W_{\text{disloc}}({\bf{\alpha}}),\delta R^T\Curl R
+R^T\Curl \delta R}\dx\nn\\
&= \int_\Omega\iprod{\DR W_{\text{mp}}(\overline{U}),
\delta R^TR\cdot R^TF}+\iprod{\DR W_{\text{disloc}}({\bf{\alpha}}),
\delta R^TR\cdot R^T\Curl R+R^T\Curl \delta R}\dx\nn\\
&= \int_\Omega\iprod{\DR W_{\text{mp}}(\overline{U})\cdot\overline{U}^T,
\delta R^TR}+\iprod{\DR W_{\text{disloc}}({\bf{\alpha}}),\delta R^TR\cdot
{\bf{\alpha}}+R^T\Curl\delta R}\dx=0\,.
\end{align}
Since $R^TR=\id$, it follows that $\delta R^TR+R^T\delta R=0$ and
$\delta R^TR=A\in\so(3)$ is arbitrary.
Therefore, Eqn.~(\ref{Euler2}) can be written as
\begin{align}
0=\int_{\Omega}\iprod{\DR W_{\text{mp}}(\overline{U})\cdot\overline{U}^T,A}
+\iprod{\DR W_{\text{disloc}}({\bf{\alpha}})\cdot{\bf{\alpha}}^T,A}
+\iprod{\DR W_{\text{disloc}}({\bf{\alpha}}),R^T\Curl(RA^T)}\;\rm{dx}
\end{align}
for all $A\in C_0^\infty(\Omega,\so(3))$. Using that $\Curl$ is a self-adjoint
operator, this is equal to
\begin{align}
0 &= \int_\Omega\iprod{\DR W_{\text{mp}}(\overline{U})\cdot
\overline{U}^T\!+\!\DR W_{\text{disloc}}({\bf{\alpha}})\,{\bf{\alpha}}^T,A}
\!+\!\iprod{\Curl(R\;\DR W_{\text{disloc}}({\bf{\alpha}})),RA^T}\dx\\
&= \int_\Omega\iprod{\DR W_{\text{mp}}(\overline{U})\cdot
\overline{U}\!+\!\DR W_{\text{disloc}}({\bf{\alpha}})\,{\bf{\alpha}}^T
\!-\!R^T\Curl(R\;\DR W_{\text{disloc}}({\bf{\alpha}})),A}\dx
\qquad\forall A\in C_0^\infty(\Omega,\so(3))\,.\nn
\end{align}
Thus, the strong form of the Euler-Lagrange equations reads
\begin{align}
\nn \text{Div}[R\;\DR W_{\text{mp}}(\overline{U})] &= 0\,,\hspace{3cm}
\text{"balance of forces"}\,,\\
\nn\skew[R^T\,\Curl(R\;\DR W_{\text{disloc}}({\bf{\alpha}}))] &=
\skew(\DR W_{\text{mp}}(\overline{U})\cdot\overline{U}^T
+\DR W_{\text{disloc}}({\bf{\alpha}})\cdot {\bf{\alpha}}^T)\,,\\
& \hspace{4cm}\text{"balance of angular momentum"}\,.
\end{align}
If $\DR W_{\text{disloc}}({\bf{\alpha}})\equiv0$ (no moment stresses,
zero characteristic length $L_c=0$) then
balance of angular momentum turns into the symmetry constraint
\begin{align}
\label{DRW}
\DR W_{\text{mp}}(\overline{U})\cdot\overline{U}^T\in\mathrm{Sym}(3)\,.
\end{align}
A complete discussion of the solutions \cite{fischle2017relaxed} to this
constraint and applications can be found in
\cite{borisov2019optimality,neff2019explicit,neff2008symmetric,
birtea2020characterization}.

\begin{figure}[h!tb]
\unitlength1cm
\begin{picture}(11.5,6.5)
\put(4.0,-0.4){\psfig{figure=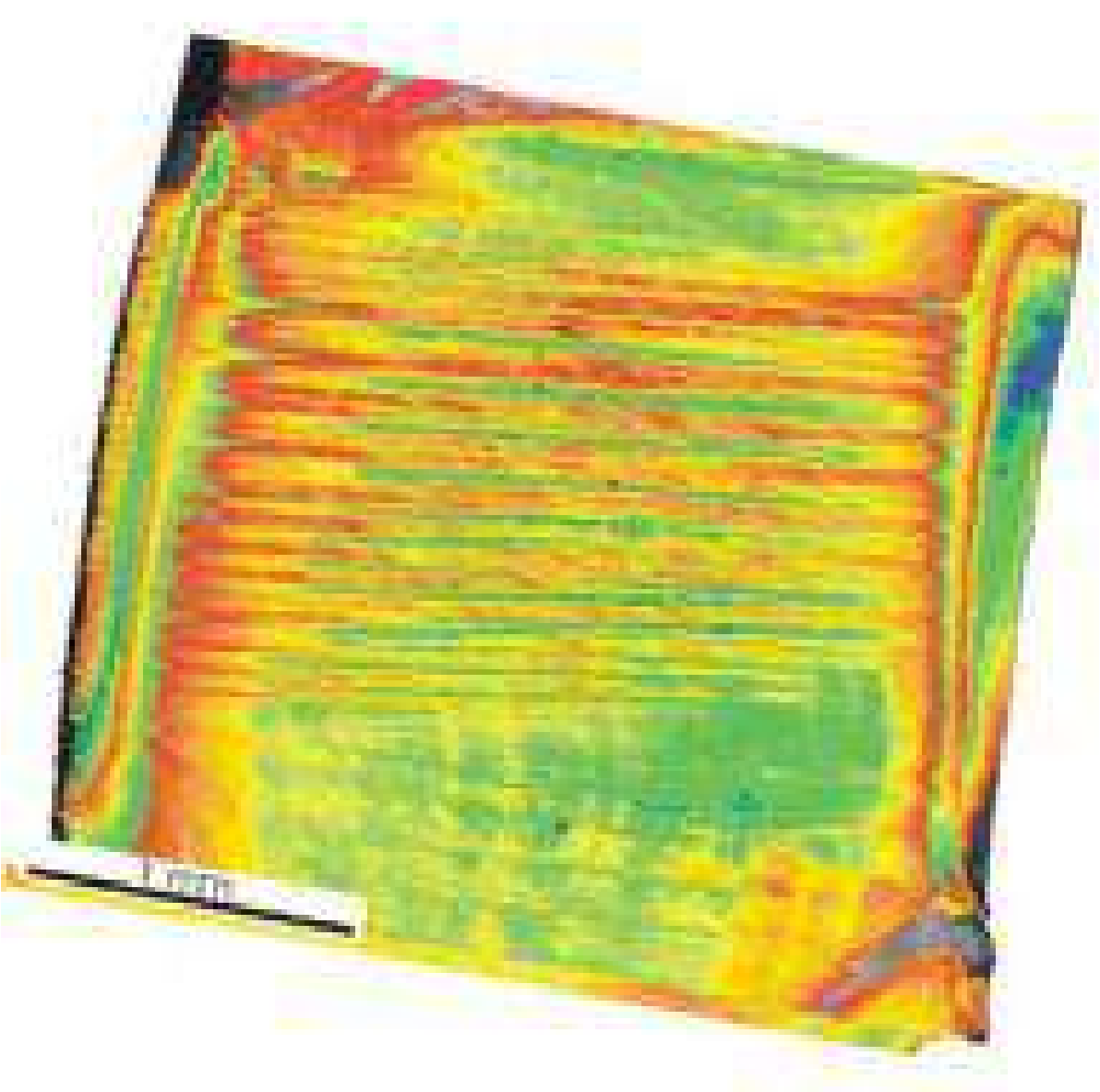,width=7.0cm}}
\end{picture}
\caption{\label{fig1}\footnotesize
A single crystal copper specimen in simple shear,
showing glide planes and micro bands. Lattice rotations do not coincide with
continuum rotations. Courtesy of D. Raabe, MPI-Eisenforschung, D\"{u}sseldorf
\cite{dmitrieva2009lamination}.}
\end{figure}

\section{ The Cosserat model in simple shear}
\begin{figure}[h!tb]
\unitlength1cm
\begin{picture}(11.5,5.0)
\put(3.5,0.0){\psfig{figure=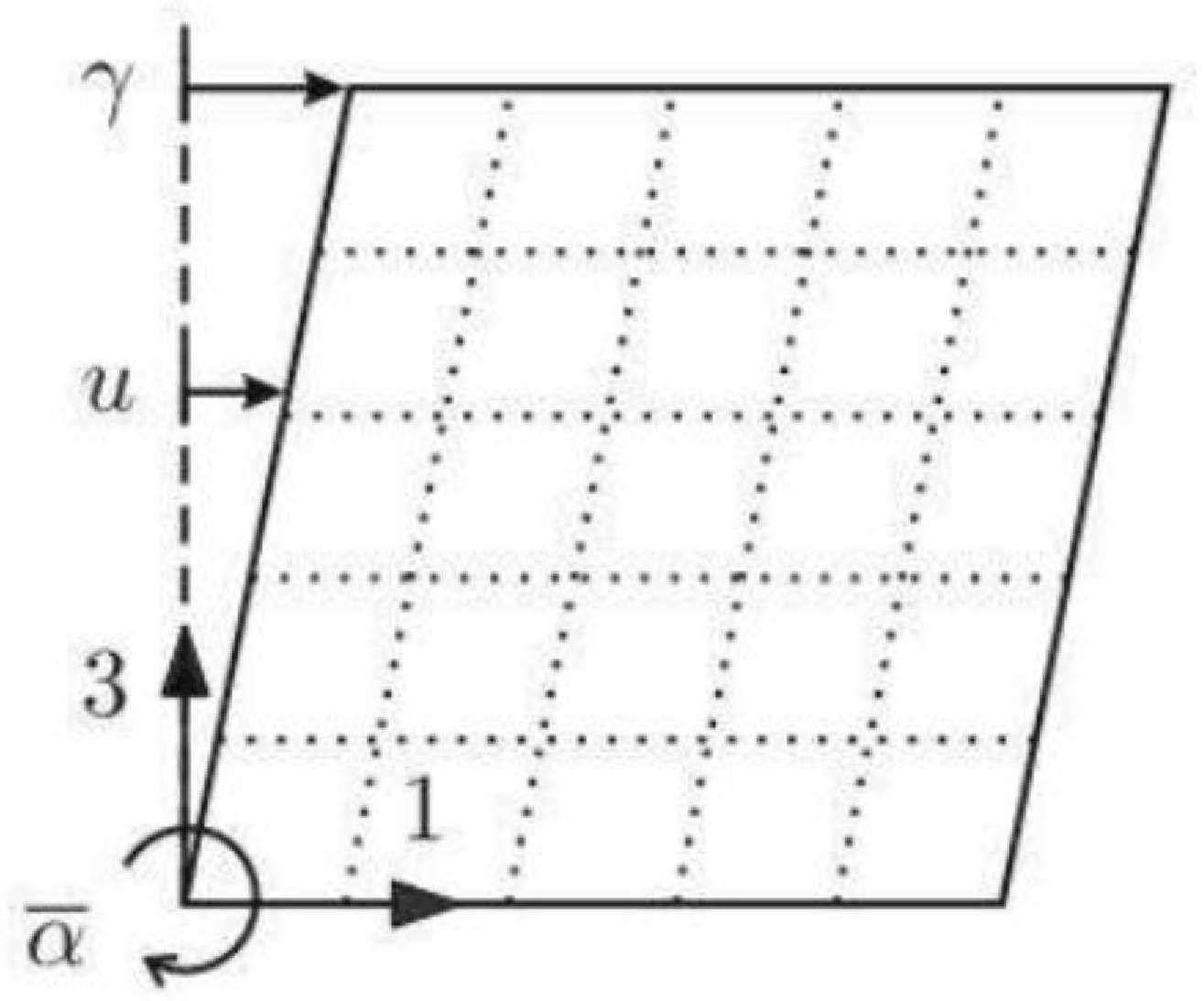,width=8.3cm}}
\end{picture}
\caption{\label{fig2}
\footnotesize The deformed state exhibits a homogeneous region in the
interior of the structure which motivates the kinematics of simple shear.}
\end{figure}

In order to elucidate the proposed nonlinear theory, notably the impact of
boundary and side conditions on the microrotations, we consider the deformation
of an infinite layer of material with unit height, fixed at the bottom and
sheared in $e_1$-direction with amount $\gamma$ at the upper face. We impose
the boundary conditions $\varphi(x_1,x_2,0)=(x_1,x_2,0)^T$,
$\varphi(x_1,x_2,1)=(x_1+\gamma,x_2,x_3)^T$, $x_1,x_2\in \R$. The parameter
$\gamma\geq 0$ is the amount of maximal shear at the upper face per unit
length. The most general deformations are of the form
\[ \varphi(x_1,x_2,x_3)=(x_1+u(x_1,x_2),x_2,x_3+v(x_1,x_3))^T, \]
see Fig~\ref{fig2}. Hence, we look for energy minimizing deformations
in the form
\begin{align}
\label{phi}
\varphi(x_1,x_2,x_3)=\begin{pmatrix}
x_1\!+\!u(x_1,x_3)\\x_2\\x_3\!+\!v(x_1,x_3)
\end{pmatrix}\,,\quad
\DR\varphi(x_1,x_2,x_3)=\begin{pmatrix}
1+u_{x_1}(x_1,x_3)&0&u_{x_3}(x_1,x_3)\\0&1&0\\
v_{x_1}(x_1,x_3)&0&1+v_{x_3}(x_1,x_3)
\end{pmatrix}\,,
\end{align}
with $u(x_1,0)=0, u(x_1,1)=\gamma$. The infinite extension in $e_1$-direction
implies that $\partial_{x_1}$ must vanish and from symmetry of the boundary
conditions at the upper and lower face, there is no reason for a displacement
in $e_3$-direction either. Hence the reduced kinematics
\begin{align}
\label{Fdef}
\varphi(x_1,x_2,x_3)=\begin{pmatrix}
x_1+u(x_3)\\x_2\\x_3
\end{pmatrix}\,,\quad
F=\DR\varphi(x_1,x_2,x_3)=\begin{pmatrix} 1&0&u'(x_3)\\0&1&0\\0&0&1
\end{pmatrix}\,,
\end{align}
with $u(0)=0$, $u(1)=\gamma$ suffices. The considered problem is therefore the
exact formulation of the simple shear in $e_1$-direction with amount $\gamma$
at the upper face of a layer of material with unit height, fixed at the bottom.

Accordingly, we assume microrotations $R\in \SO(3)$ in the form
\begin{align}
\label{R}
R(x_1,x_2,x_3)=\begin{pmatrix}
\cos\alpha(x_3)&0&\sin\alpha(x_3)\\0&1&0\\
-\sin\alpha(x_3)&0&\cos\alpha(x_3) \end{pmatrix}\,,
\end{align}
having fixed axis of rotation $e_2$. Therefore,
\begin{align}
\Curl R=\begin{pmatrix}
0 & -\sin\alpha(x_3)\alpha'(x_3) & 0\\ 0 & 0 &0\\
0 & -\cos\alpha(x_3)\alpha'(x_3) & 0 \end{pmatrix}\,.
\end{align}
In the following, we denote $x_3$ by $x$. It holds
$|\overline{R}^T\Curl\overline{R}|^2=|\Curl\overline{R}|^2
=|\overline{\alpha}'|^2$.

Inserting the ansatz (\ref{phi}) and (\ref{R}) leaves us with the energy
\begin{align}
E(u,\alpha) \;=\;& \int_0^1 W(u',\alpha,\alpha') \dx\nn\\
:=\; & \mu\!\!\io\!2L_c^2|\alpha'|^2\!+\!2(\ca\!-\!1)^2\!+\!
\frac{1\!+\!\ssa}{2}|u'|^2\!+\!2(\ca\!-\!1\!)\sa u'\dx\nn\\
&\quad +\frac{\mu_c}{2}\io\cca\big(2\ta-u'\big)^2\dx.
\end{align}

\section{A one-dimensional simple shear problem}
\label{secshear}
In this article we are now concerned with minimizers of the
{\it mechanical energy functional}
\begin{align}
E(u,\alpha) \;&=\;\; \mu\!\!\io\!2L_c^2|\alpha'|^2\!+\!2(\ca\!-\!1)^2\!+\!
\frac{1\!+\!\ssa}{2}|u'|^2\!+\!2(\ca\!-\!1\!)\sa u'\dx\nn\\
\label{Edef}
&\; \qquad +\frac{\mu_c}{2}\io\cca\big(2\ta-u'\big)^2\dx.
\end{align}

Introducing the third-order expansion $\ca\sim1-\frac{\alpha^2}{2}$,
$\sa\sim\alpha-\frac{\alpha^3}{6}$ and dropping the higher order terms
leads to the {\it reduced mechanical energy functional}
\begin{align}
\Ered(u,\alpha)
\;:=\;\;& \mu\io\Big(2L_c^2|\alpha'|^2+\frac{1+\alpha^2}{2}
|u'|^2+\frac{\alpha^4}{2}-\alpha^3u'\Big)\dx\nn\\
\label{Wred}
& +2\mu_c\int_0^1\Big(\frac{u'}{2}-\alpha\Big)
\Big(\frac{u'}{2}-\alpha-\frac{\alpha^2}{6}(3u'-2\alpha)\Big)\dx.
\end{align}
For the definitions (\ref{Edef}) and (\ref{Wred}), see also
\cite{neff2009simple}[Eqn.(3.11), Eqn.(3.22)].
We first provide alternative representations of these two functionals as this
allows to simplify the Euler-Lagrange equations and will give insights into the
minimizers later.

\begin{Lemma}
\label{lem1}
The functionals $E$, $\Ered$ defined in (\ref{Edef}), (\ref{Wred})
can alternatively be written as
\begin{align}
\label{Wdef}
E(u,\alpha) \;=\;\; & \io W(u',\alpha,\alpha')\dx\\
\;=\;\; &\frac{\mu}{2}\int_0^14L_c^2|\alpha'|^2+|u'|^2
+\Big(\sa u'-4\shalf\!\Big)^2\dx\nn\\
\label{Edef2}
& +\frac{\mu_c}{2}\int_0^1\Big(\ca u'-2\sa\Big)^2\dx,
\end{align}
\begin{align}
\Ered(u,\alpha) \;=\;\;
& \frac{\mu}{2}\int_0^14L_c^2|\alpha'|^2+|u'|^2
+[\alpha(\alpha\!-\!u')]^2\dx\nn\\
\label{Wred2}
& +\frac{\mu_c}{2}\int_0^1(1-\alpha^2)|u'|^2+\Big(\frac{8\alpha^2}{3}
-4\Big)\alpha u'+4\alpha^2-\frac{4\alpha^4}{3}\dx.
\end{align}
\end{Lemma}

\noindent{\bf Proof.} Both identities follow from straightforward elementary
rearrangements. Trigonometric addition formulas imply
$\cos(2\alpha)=\cca-\ssa=1-2\ssa$, hence
\begin{equation}
\label{s0}
\ca-1=-2\shalf.
\end{equation}
Plugging (\ref{s0}) into (\ref{Edef}) yields
\begin{eqnarray}
&& \hspace*{-50pt}
\nonumber 2(\ca-1)^2+\frac{1+\ssa}{2}|u'|^2+2(\ca-1)\sa\,u'\\
\nonumber&=& 8\sin^4\!\big(\frac{\alpha}{2}\big)+\frac{|u'|^2}{2}
+\frac{\ssa|u'|^2}{2}-4u'\shalf\sa\\
&=& \frac12\Big(\sa u'-4\shalf\Big)^2+\frac{|u'|^2}{2}.
\end{eqnarray}
Similarly, for the $\mu_c$-integral in (\ref{Edef}),
\begin{eqnarray}
\nonumber&& \hspace*{-40pt}
\cca(2\ta-u')^2=\cca\big(4\tan^2(\alpha)-4\ta u'+|u'|^2\big)\\
&=& 4\ssa-4\sa\ca u'+\cca|u'|^2=\Big(\ca u'-2\sa\Big)^2.
\end{eqnarray}
This proves (\ref{Edef2}). The proof of (\ref{Wred2}) is immediate from
\[ \frac{1+\alpha^2}{2}|u'|^2+\frac{\alpha^4}{2}-\alpha^3u'=
\frac{|u'|^2}{2}+\frac{\alpha^2}{2}\Big(|u'|^2-2\alpha u'+\alpha^2\Big)\,, \]
and a re-ordering of the $\mu_c$-integral. \qed
\vspace*{4mm}

\begin{Remark}
\label{rem1}
The $\mu_c$-part of $\Ered$ cannot be written as a complete quadratic form since
some higher order terms have been dropped. However, by keeping the terms
$\frac{\mu_c}{2}\Big(\frac{\alpha^4}{4}|u'|^2-\frac13\alpha^5u'
+\frac19\alpha^6\Big)$ of the Taylor expansion of $W$, we find
\begin{align}
\label{Wreddef}
\Ered(u,\alpha) \;&=\;\; \io\Wred(u',\alpha,\alpha')\dx\\
\label{Wred3}
&=\; \frac{\mu}{2}\int_0^1\!4L_c^2|\alpha'|^2+|u'|^2
+[\alpha(\alpha\!-\!u')]^2\dx+\frac{\mu_c}{2}\!\int_0^1\!\Big(
\frac{2\!-\!\alpha^2}{2}u'-\frac{6\alpha\!-\!\alpha^3}{3}\Big)^2\dx.
\end{align}
From now on we will use this representation (\ref{Wred3}) instead of
(\ref{Wred2}).
\end{Remark}

\vspace*{3mm}
We consider the minimization problem
\begin{equation}
\label{Wmin}
E(u,\alpha)\to\min
\end{equation}
subject either to the {\it consistent coupling conditions}
\cite{d2021consistent} (derived from $\mathrm{skew}(R^TF)_{|\{0,1\}}=0$)
\begin{equation}
\label{cc}
\hspace*{-50pt} u(0)=0,\,u(1)=\gamma,\qquad\qquad
u'(0)=2\tan(\alpha(0)),\quad u'(1)=2\tan(\alpha(1))
\end{equation}
or subject to the {\it Dirichlet boundary conditions}
(for prescribed microrotation angle $\alpha_D\in\R$ at the upper and
lower faces)
\begin{equation}
\label{BC}
u(0)=0,\;u(1)=\gamma,\qquad\qquad \alpha(0)=\alpha(1)=\alpha_D.
\end{equation}
In order to single out solutions, we may impose the further
{\it symmetry constraint}
\begin{equation}
\label{per}
u'(1)=u'(0).
\end{equation}

The Euler-Lagrange equations related to (\ref{Wmin}), replacing
\cite{neff2009simple}[Eqn.~(3.14)], read
\begin{align}
\label{S1}
[2\mu\!+\!(\mu_c\!\!-\!\!\mu)\cos^2(\alpha)]u'' \;=\;\; & 2(\mu_c\!\!-\!\!\mu)
\Big(\sa\ca u'\!+\!\cca\!-\!\ssa\Big)\alpha'\!+\!2\mu\ca\alpha',\\
\label{S2}
4\mu L_c^2\alpha'' \;=\;\; &
\big(\!\ca u'\!-\!2\sa\big)\big[(\mu\!-\!\mu_c)\big(
\sa u'\!-\!4\shalf\!\big)\!-\!2\mu_c\big]
\end{align}
for $x=x_3\in\Omega:=(0,1)$
subject to either (\ref{cc}), (\ref{per}) or
(\ref{BC}), (\ref{per}). The equations (\ref{S1}), (\ref{S2}) constitute
balance of force and angular momentum, respectively.
Eqn.~(\ref{S2}) is the form that Eqn.~(\ref{DRW}) takes for the
ansatz made here.

Eqn.~(\ref{S1}) is equivalent to $\ddx\tau(u,\alpha)=0$ in $\Omega$ for the
{\it force stress tensor} $\tau=\tau(u,\alpha)$ defined by
\begin{equation}
\label{taudef}
\tau:=\DR_{u'}W(u',\alpha,\alpha')=\mu\!\Big[u'+\sa\!\Big(\sa u'-4\shalf
\Big)\Big]+\mu_c\cos(\alpha)\!\Big[\ca u'\!-\!2\sa\Big].
\end{equation}
The Euler-Lagrange equations related to $\Ered(u,\alpha)\to\min$ lead to the
{\it reduced system}
\begin{align}
-\Big[\mu(1\!+\!\alpha^2)+\mu_c\Big(1\!-\!\alpha^2
\underline{+\frac14\alpha^4}\Big)\Big]u'' \;=\; & (4\mu_c\!-\!3\mu)
\alpha^2\alpha'+2(\mu\!-\!\mu_c)\alpha\alpha'\,u'-2\mu_c\alpha'\nn\\
\label{red1}
&\qquad \underline{+\mu_c\alpha^3\Big(u'-\frac56\alpha\Big)\alpha'},\\
\label{red2}
\mu\Big(\!\!-L_c^2\alpha''+\frac12\alpha^3-\frac34\alpha^2 u'
+\frac14\alpha|u'|^2\!\Big) &
+\mu_c\Big(\!\!-\frac14\alpha|u'|^2+\alpha^2u'
-\frac23\alpha^3+\alpha-\frac12u'\!\Big)\nn\\
&\qquad +\underline{\mu_c\Big(\frac18\alpha^3|u'|^2-\frac{5}{24}\alpha^4u'
+\frac{1}{12}\alpha^5\Big)}=0
\end{align}
for $x=x_3\in\Omega:=(0,1)$ subject either to the reduced consistent coupling
boundary conditions
\begin{equation}
\label{ccred}
\hspace*{-30pt}
u(0)=0,\;u(1)=\gamma,\qquad\qquad 2\alpha(0)=u'(0),\;2\alpha(1)=u'(1)
\end{equation}
and (\ref{per}) or subject to (\ref{BC}) and (\ref{per}).
Underlined in (\ref{red1}), (\ref{red2}) are those higher order terms that are
only present if $\Ered$ is defined by (\ref{Wred3}) instead of (\ref{Wred})
or (\ref{Wred2}).

\begin{figure}[h!t]
\unitlength1cm
\begin{picture}(11.5,3.0)
\put(0.1,-0.30){\psfig{figure=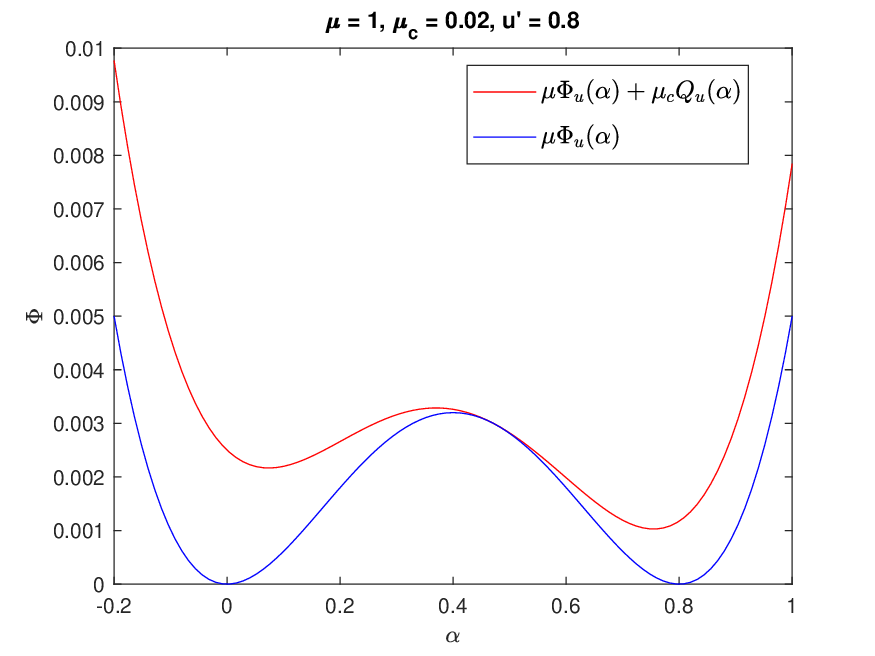,width=5.2cm}}
\put(5.5,-0.30){\psfig{figure=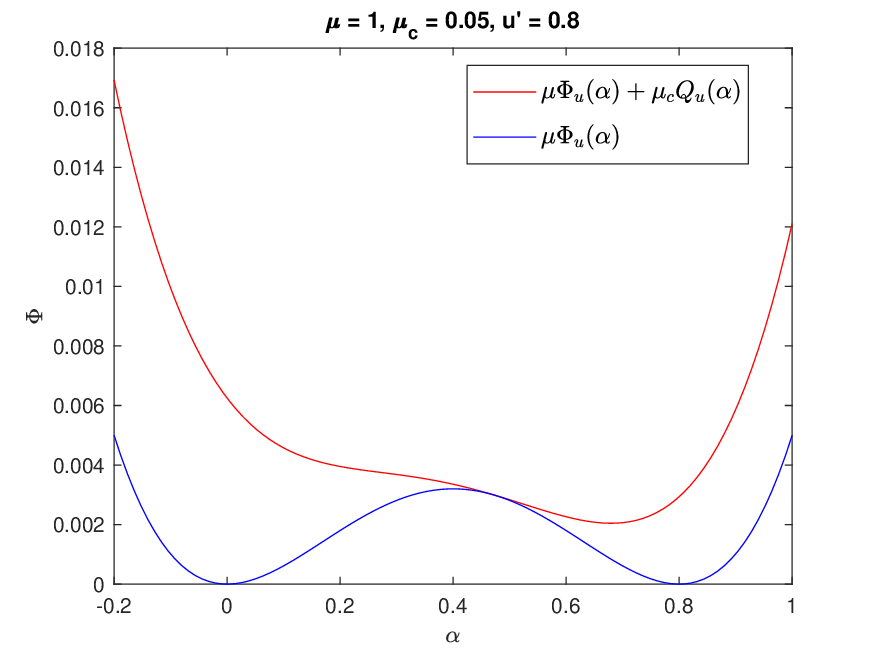,width=5.2cm}}
\put(11.0,-0.30){\psfig{figure=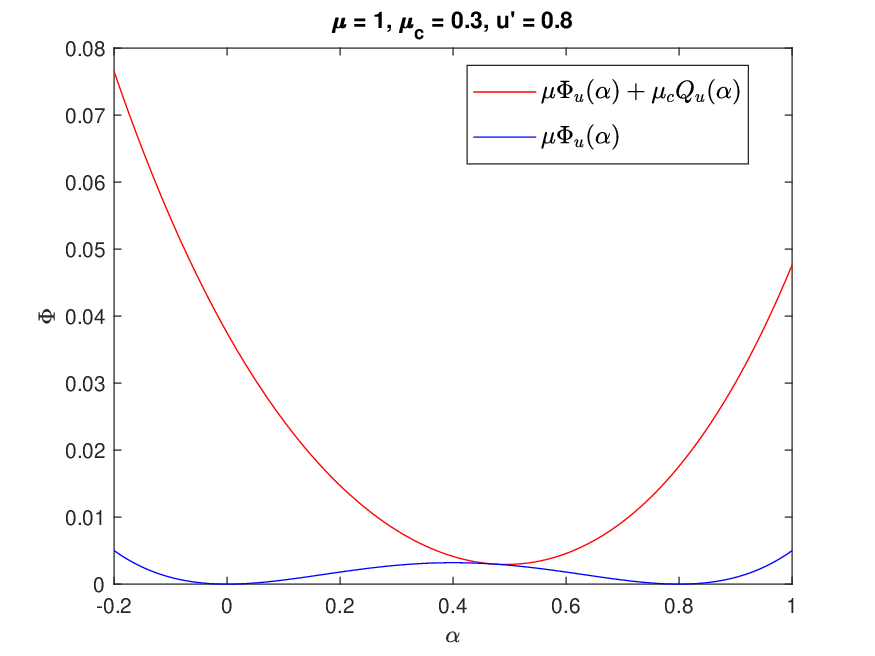,width=5.2cm}}
\end{picture}
\caption{\label{fig3}
The potential $\mu\,\Phi_u(\alpha)+\mu_c\,Q_u(\alpha)$ (in red) compared to the
double-well potential $\mu\,\Phi_u(\alpha)$ (in blue) for
$\alpha\in[-0.2,1]$, $u'=0.8$ and $\mu=1$.
Left: $\mu_c=0.02$. Center: $\mu_c=0.05$. Right: $\mu_c=0.3$.}
\end{figure}

\noindent
Eqn.~(\ref{red1}) is equivalent to $\ddx\taured(u,\alpha)=0$ in $\Omega$ for
the {\it reduced force stress tensor} $\taured=\taured(u,\alpha)$
defined by
\begin{equation}
\label{taureddef}
\taured:=\DR_{u'}W(u',\alpha,\alpha')=\mu\Big[u'
-\alpha^2(\alpha\!-\!u')\Big]+\mu_c\Big(\frac{2\!-\!\alpha^2}{2}\Big)\Big[
\frac{2\!-\!\alpha^2}{2}u'-\frac{6\alpha\!-\!\alpha^3}{3}\Big].
\end{equation}
We point out that Eqn.~(\ref{red2}) can equally be written as
\begin{equation}
\label{AC}
\mu L_c^2\alpha'' = \mu\,\partial_\alpha\Phi_u(\alpha)+\mu_c\,\partial_\alpha
Q_u(\alpha)
\end{equation}
with the potentials
\begin{equation}
\Phi_u(\alpha) := \frac18\alpha^2(\alpha-u')^2,\qquad\qquad
Q_u(\alpha) := \frac18\Big(\frac{2-\alpha^2}{2}u'-\frac{6\alpha-\alpha^3}{3}
\Big)^2.
\end{equation}
For small values $\mu_c\ge0$, the double-well potential $\Phi_u$
is dominant and Eqn.~(\ref{AC}) corresponds to a stationary Allen-Cahn
equation\footnote{Allen-Cahn equation: $\pat u=\varepsilon\Delta u-\psi'(u)$
with a double-well potential $\psi$.}. The larger the Cosserat couple modulus
$\mu_c$, the stronger the influence of the quadratic potential
$Q_u$. There is a bifurcation and a critical value $\mu_c^{\rm crit}$, such
that for $\mu_c\ge\mu_c^{\rm crit}$, the right hand side of Eqn. (\ref{AC})
has only one minimizer, see Figs.~\ref{fig3}.

\begin{figure}[pht]
\unitlength1cm
\begin{picture}(11.5,3.45)
\put(5.5,-0.40){\psfig{figure=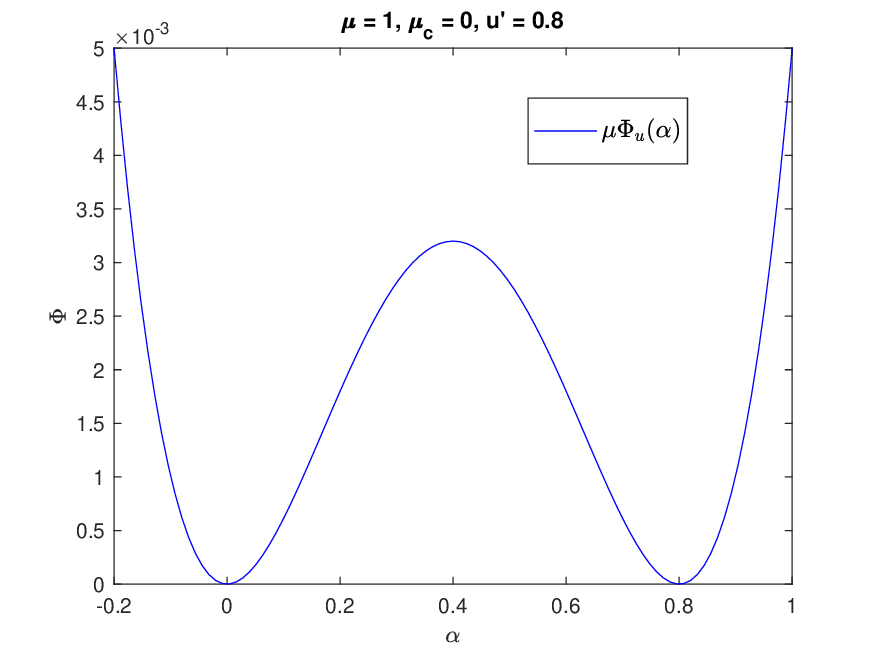,width=5.2cm}}
\end{picture}
\caption{\label{fig3b}
The double-well potential $\mu\Phi_u(\alpha)$ for $\alpha\in[-0.2,1]$,
$u'=0.8$ and $\mu=1$.}
\end{figure}

\vspace*{4mm}
Related to different boundary conditions we introduce the reflexive
Banach spaces
\begin{align}
\nonumber\cXD &:= \Big\{(u,\alpha)\!\in\!(W^{1,2}(\Omega;\,\R))^2\,\Big|
\,u(0)\!=\!0,\,u(1)\!=\!\gamma, \quad
\alpha(0)\!=\!\alpha(1)\!=\!\alpha_D\Big\},\\
\XD &:= \Big\{(u,\alpha)\!\in\!W^{1,2}(\Omega;\,\R)\!\times\!
L^4(\Omega;\,\R)\,\Big|\,u(0)\!=\!0,\,u(1)\!=\!\gamma,\quad
\alpha(0)\!=\!\alpha(1)\!=\!\alpha_D\Big\}
\end{align}
and correspondingly for the consistent coupling conditions
\begin{align}
\begin{split}
\cXCred &:= \Big\{(u,\alpha)\!\in\!(W^{1,2}(\Omega;\,\R))^2\,\Big|\,
u(0)\!=\!0,\,u(1)\!=\!\gamma,\,2\alpha(0)\!=\!u'(0),
\,2\alpha(1)\!=\!u'(1)\Big\},\\
\cX^C &:= \Big\{(u,\alpha)\!\in\!(W^{1,2}(\Omega;\,\R))^2\,\Big|\,
u(0)\!=\!0,\,u(1)\!=\!\gamma, \,2\tan\alpha(0)\!=\!u'(0),\,
2\tan\alpha(1)\!=\!u'(1)\Big\},\\
\XCred &:= \Big\{(u,\alpha)\!\in\!W^{1,2}(\Omega;\,\R)\!\times
\!L^4(\Omega;\,\R)\,\Big|\,u(0)\!=\!0,\,u(1)\!=\!\gamma,\,
2\alpha(0)\!=\!u'(0),\,2\alpha(1)\!=\!u'(1)\Big\},
\end{split}
\end{align}
\vspace*{-6mm}
\begin{align*}
X^C &:= \Big\{(u,\alpha)\!\in\!W^{1,2}(\Omega;\,\R)\!\!\times\!
L^4(\Omega;\,\R)\,
\Big|\,u(0)\!=\!0,\,u(1)\!=\!\gamma,\,2\tan\alpha(0)\!=\!u'(0),
2\tan\alpha(1)\!=\!u'(1)\Big\}.
\end{align*}

\vspace*{3mm}
\begin{Lemma}
\label{lem2}
Let $L_c>0$. Then, for any $\mu>0$, $\mu_c\ge0$, $\Ered$ defined by
(\ref{Wred3}) possesses a minimizer $(u,\alpha)$ in $\cXCred$, $\cX^D$ and
$E$ possesses a minimizer in $\cX^C$, $\cX^D$. For $L_c=0$, $\Ered$ has
minimizers in $\XCred$, $X^D$ and $E$ has minimizers in $X^C$, $X^D$.
\end{Lemma}

\noindent{\bf Proof.} (i) Let $L_c>0$. We first consider the case of Dirichlet
boundary conditions for $u$ and $\alpha$. Rewriting (\ref{Wred3}) as a
functional defined on $\cX_0:=(W^{1,2}_0((0,1);\;\R)^2$, we find with Young's
inequality for $(u,\alpha)\in\cX_0$
\begin{align}
\Ered(u+\gamma,\alpha+\alpha_D) &=\; \frac{\mu}{2}\int_0^14L_c^2
|\alpha'|^2+(u'\!+\!\gamma)^2+\big[(\alpha\!+\!\alpha_D)(\alpha\!+\!\alpha_D
-(u'\!+\!\gamma))\big]^2\dx\nn\\
&\quad\quad +\frac{\mu_c}{2}\int_0^1\Big(\frac{2-(\alpha\!+\!\alpha_D)^2}{2}
(u'\!+\!\gamma)-\frac{6(\alpha\!+\!\alpha_D)-(\alpha\!+\!\alpha_D)^3}{3}
\Big)^2\dx\nn\\
\label{Wcoerc}
&\ge\; \frac{\mu}{2}\int_0^14L_c^2|\alpha'|^2+(u'\!+\!\gamma)^2\dx\\
&  =\mu\int_0^12L_c^2|\alpha'|^2+\frac12|u'|^2+\big(\frac{u'}{2}\big)
(2\gamma)+\frac{\gamma^2}{2}\dx\nn\\
&\ge\; \mu\int_0^12L_c^2|\alpha'|^2+\frac12|u'|^2-\frac{1}{4}|u'|^2
-4\gamma^2+\frac{\gamma^2}{2}\dx\nn\\
& =\mu\int_0^12L_c|\alpha'|^2+\frac14|u'|^2-\frac{7\gamma^2}{2}\dx.\nn
\end{align}
With the Poincar{\'e} inequality valid on $\cX_0$ and the Banach-Alaoglu theorem
this demonstrates that the level sets
\[ \Big\{(u,\alpha)\in\cX^D\;\Big|\;\Ered(u,\alpha)\le C\Big\} \]
for constants $C>0$ are sequentially weakly precompact proving the coercivity
of $\Ered$, see, e.g. \cite{Struwe}.

We observe that the integrand $\Wred(u',\alpha,\alpha')$ as defined in
(\ref{Wreddef}) is a Carath{\'e}odory function and strictly convex both in
$u'$ and $\alpha'$.
Despite the dependence of $\Wred$ on $\alpha$, the proof of weak lower
semicontinuity of $\Ered$ in $\cX^D$ can thus be carried out in the spirit of
the well-known Tonelli-Serrin theorem, see \cite[Section~3.2.6]{Dac08}
for details.
Alternatively, the weak lower semicontinuity can be derived from the more
general result in \cite{AF84} based on gradient Young measures.
By the direct method in the calculus of variations, the coercivity and weak
lower semicontinuity of $\Ered$ yield the existence of a minimizer
$(u,\alpha)\in\cX^D$. The proof of minimizers of $\Ered$ in $\cXCred$ is
similar.

Now let us consider the case of Dirichlet boundary conditions for $E$.
Proceeding as above and estimating the quadratic terms from below by $0$, we
find for $(u,\alpha)\in\cX_0$
\[ E(u+\gamma,\alpha+\alpha_D)\ge\frac{\mu}{2}\int_0^14L_c^2|\alpha'|^2
+(u'+\gamma)^2\dx. \]
This coincides with (\ref{Wcoerc}). From there, with the Poincar{\'e}
inequality the coercivity of $E$ in $\cX^D$ can be shown as above. The lower
semicontinuity of $E$ is again a consequence of strict convexity and the
Tonelli-Serrin theorem. The proof of minimizers of $E$ in $\cX^C$ is similar.

(ii) Let $L_c=0$.
The coercivity of $\Ered$ w.r.t. weak-convergence can be shown similar to (i).
In contrast, $E$ is $2\pi$-periodic in $\alpha$ which replaces the coercivity
in $\alpha$.
Eqns.~(\ref{Wred3}), (\ref{Edef2}) imply the lower semicontinuity of $\Ered$,
$E$. With the direct method, the existence of minimizers of $\Ered$ in $X^D$,
$\XCred$ and of $E$ in $X^D$, $X^C$ for any $\mu>0$, $\mu_c\ge0$ follows. \qed

\begin{figure}[h!bt]
\unitlength1cm
\begin{picture}(11.5,3.6)
\put(5.5,-0.20){\psfig{figure=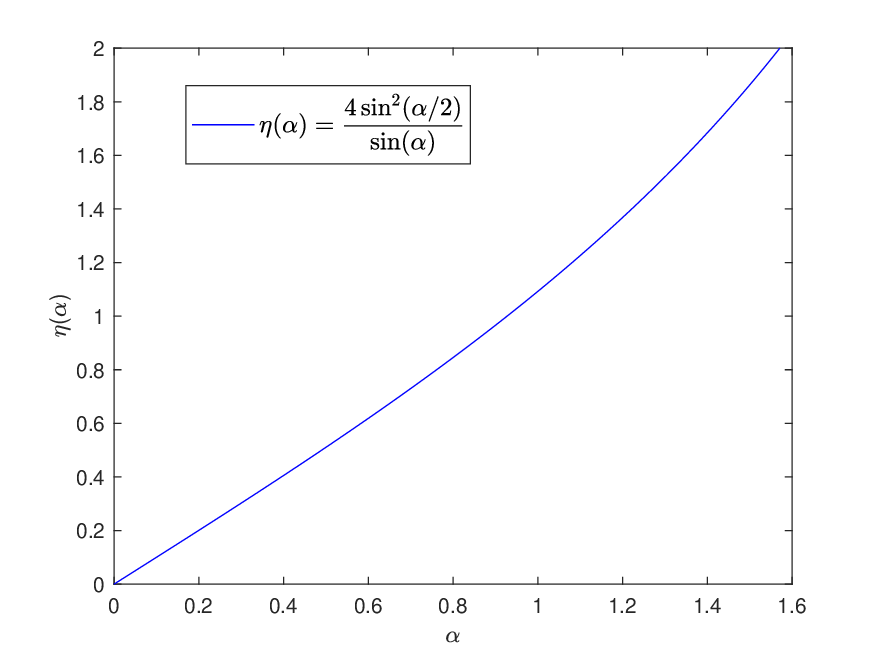,width=5.3cm}}
\end{picture}
\caption{\label{fig4}
Plot of the function $\eta(\alpha)=\frac{4\shalf}{\sa}$ for
$\alpha\in\big[0,\frac\pi2\big]$.}
\end{figure}

\newpage
For the detailed discussion of the Euler-Lagrange equations related to problem
(\ref{Wmin}) we introduce the function
\begin{equation}
\label{etadef}
\eta(\alpha):=\frac{4\shalf}{\sa}.
\end{equation}
It holds $\lim_{\alpha\to0}\eta(\alpha)=0$ and $\eta(\alpha)$ is invertible,
monotone increasing and close to a linear function, cf. Fig.~\ref{fig4}.
 
In the following discussion, we first ignore the boundary conditions on
$u$ and $\alpha$ and study sufficiently regular solutions
$u\in W^{2,2}(\Omega)$ and $\alpha\in W^{1,2}(\Omega)$.
Clearly, not every solution is in this class, c.f. Section~\ref{secnum}.

\begin{Lemma}
\label{lem3}
Let $L_c=0$ and $\mu>0$, $\mu_c\ge0$. Then every solution
\[ (u,\alpha)\in W^{2,2}(\Omega)\times W^{1,2}(\Omega) \]
of the Euler-Lagrange equations (\ref{S1}), (\ref{S2}) is continuous in
$\overline{\Omega}$ and satisfies
\begin{equation}
\label{uzero}
u''(x)=\alpha'(x)=0 \mbox{ for a.e. }x\in\Omega.
\end{equation}
\end{Lemma}

\noindent{\bf Proof.} Due to the Sobolev embedding
$W^{1,2}(\Omega)\hookrightarrow C^0(\overline{\Omega})$, both $\alpha$ and $u$
are continuous.

(i) For $\mu=\mu_c$, the Euler-Lagrange equations (\ref{S1}), (\ref{S2})
simplify to
\begin{align}
\label{EL1a}
u'' \;&=\; \ca\alpha',\\
\label{EL1b}
0 \;&=\; \ca u'-2\sa
\end{align}
to be satisfied pointwise for $x\in\Omega$.

Eqn.~(\ref{EL1b}) is equivalent to $u'=2\ta$ and taking the derivative
yields $u''=\frac{2\alpha'}{\cca}$. With (\ref{EL1a}) we find
\begin{equation}
\label{x1}
\frac{2\alpha'}{\cca}=\ca\alpha'.
\end{equation}
Since $2=\cos^3(\alpha)$ has no solution, Eqn.~(\ref{x1}) implies $\alpha'=0$,
and so with (\ref{EL1a}) $u''=0$.
\newpage

\vspace*{5mm}
\noindent(ii) Let $\mu>0$, $\mu_c\ge0$ with $\mu\not=\mu_c$. Considering the
Euler-Lagrange equation (\ref{S2}) with $L_c=0$, two cases may occur.

\noindent\underline{Case~1}: $\alpha=\arctan\!\big(\frac{u'}{2}\big)$.

This implies $u'=2\ta$ such that $u''=\frac{2\alpha'}{\cca}$.
Inserting this identity in (\ref{S1}) yields
\begin{equation}
\label{s2}
(\mu_c\!-\!\mu)\frac{\big[\sa\ca u'\!+\!\cca\!-\!\ssa\big]\alpha'}
{2\mu+(\mu_c\!-\!\mu)\cca}
+\frac{\mu\ca\alpha'}{2\mu\!+\!(\mu_c\!-\!\mu)\cca}=\frac{\alpha'}{\cca}.
\end{equation}
[The denominator in (\ref{s2}) is positive, e.g.
$2\mu+(\mu_c\!-\!\mu)\cos^2(\alpha)\ge2\mu>0$ for $\mu_c\ge\mu$ and
$2\mu+(\mu_c\!-\!\mu)\cos^2(\alpha)\ge2\mu-(\mu_c\!-\!\mu)>0$ for $\mu_c<\mu$.]

Using $u'=2\ta$ on the right gives after simplifications
\begin{equation}
\label{no1}
\frac{\big[\mu_c-\mu+\mu\ca\big]\alpha'}{2\mu+(\mu_c\!-\!\mu)\cca}
=\frac{\alpha'}{\cca}
\end{equation}
which is equivalent to
\[ \mu\cos^3(\alpha)\alpha'=2\mu\alpha'. \]
Since $\mu>0$, this proves as in (i) that $\alpha'=0$ which with (\ref{S1})
results in $u''=0$.

\vspace*{2mm}
\noindent\underline{Case~2}:
$\sa u'=4\shalf+\frac{2\mu_c}{\mu-\mu_c}$.

For $\mu_c=0$, one solution is $\alpha\equiv2\pi k$ for some $k\in\Z$.
Then $\alpha'=0$ and Eqn.~(\ref{S1}) shows $u''=0$. Alternatively, for
$\alpha\not\equiv2\pi k$, we have $\sin(\alpha)\not=0$ such that
\begin{equation}
\label{uder}
u'=\frac{2\mu_c}{(\mu-\mu_c)\sa}+\frac{4\shalf}{\sa}.
\end{equation}
Taking the derivative yields
\[ u''=-\frac{2\mu_c\ca\alpha'}{(\mu\!-\!\mu_c)\ssa}
+\frac{\big[4\sin\!\big(\frac{\alpha}{2}\big)\cos\!\big(
\frac{\alpha}{2}\big)\sa-4\shalf\ca\big]\alpha'}{\ssa}. \]
We use $4\sin\!\big(\frac{\alpha}{2}\big)\cos\!\big(\frac{\alpha}{2}
\big)=2\sa$ and $-4\shalf=-2+2\ca$.
With the Euler-Lagrange equation (\ref{S1}) this leads to
\begin{align*}
& \hspace*{-30pt}
-\frac{2\mu_c\ca\alpha'}{(\mu\!-\!\mu_c)\ssa}
+\frac{2\big[\ssa-\ca+\cca\big]\alpha'}{\ssa}\\
& \;=\; \frac{\big[2(\mu_c\!-\!\mu)\big(\sa\ca u'+\cca-\ssa\big)+2\mu\ca\big]
\alpha'}{2\mu+(\mu_c\!-\!\mu)\cca}
\end{align*}
and simplifies to
\begin{align*}
& \hspace*{-40pt}
\frac{\mu_c\ca\alpha'}{(\mu_c\!-\!\mu)(1-\cca)}
+\frac{\big[1-\ca\big]\alpha'}{1-\cca}\\
& \;=\; \frac{\big[(\mu_c\!-\!\mu)\big(\sa\ca u'+\cca-\ssa\big)+\mu\ca\big]
\alpha'}{2\mu+(\mu_c\!-\!\mu)\cca}.
\end{align*}
On the right we plug in the expression (\ref{uder}) for $u'$. So we obtain
\begin{align*}
& \frac{\big[(\mu_c\!-\!\mu)\big(1-\ca\big)+\mu_c\ca\big]\alpha'}
{(\mu_c\!-\!\mu)\big(1-\cca\big)}\\
& \; = \; \frac{\big[-2\mu_c\ca+\mu\ca+(\mu_c\!-\!\mu)\big(2\ca-2\cca+\cca
-\ssa\big)\big]\alpha'}{2\mu+(\mu_c\!-\!\mu)\cca}.
\end{align*}
After simplifications, this is equivalent to
\begin{equation}
\label{c2simp}
\frac{\big[(\mu_c\!-\!\mu)+\mu\ca\big]\alpha'}{(\mu_c\!-\!\mu)\big(1-\cca
\big)}=\frac{\big[-\mu\ca-(\mu_c\!-\!\mu)\big]\alpha'}{2\mu+(\mu_c\!-\!\mu)
\cca}.
\end{equation}
Expanding this, we find
\begin{align}
& \big[2\mu(\mu_c\!-\!\mu)+2\mu^2\ca+(\mu_c\!-\!\mu)^2\cca+\mu(\mu_c\!-\!\mu)
\cos^3(\alpha)\big]\alpha'\nn\\
\label{c2simp2}
& \; = \; \big[-\mu(\mu_c\!-\!\mu)\ca-(\mu_c\!-\!\mu)^2+\mu(\mu_c\!-\!\mu)
\cos^3(\alpha)+(\mu_c\!-\!\mu)^2\cca\big]\alpha'.
\end{align}
This simplifies further to
\begin{equation}
\label{c2simp3}
[2\mu^2+\mu(\mu_c\!-\!\mu)\big]\ca\alpha'=\big[2\mu(\mu\!-\!\mu_c)
-(\mu_c\!-\!\mu)^2\big]\alpha'
\end{equation}
and eventually
\[ \mu(\mu\!+\!\mu_c)\ca\alpha'=(\mu\!+\!\mu_c)(\mu\!-\!\mu_c)\alpha'. \]
Hence either $\alpha'=0$ or
\begin{equation}
\label{rend}
\alpha\equiv\alpha_4:=\arccos\Big(\frac{\mu-\mu_c}{\mu}\Big)
\end{equation}
provided $\big|\frac{\mu-\mu_c}{\mu}\big|\le1$.
In both cases we have $\alpha'=0$ which shows with (\ref{S1}) that $u''=0$. \qed

\vspace*{2mm}
\begin{Remark}
\label{rem2}
The consistent coupling condition is made such that the homogeneous deformation
$\uhom(x)=\gamma x$ remains a solution of the boundary value problem.
\end{Remark}

\newpage
\noindent For given $\gamma>0$ we introduce the constants
\begin{equation}
\label{mindef}
\cblue\alpha_1^-:=\arctan\!\Big(\!\frac{\gamma\mu-f}{2\mu+\frac{\gamma}{2}f}\!
\Big),\quad \alpha_1^+:=\arctan\!\Big(\!\frac{\gamma\mu+f}{2\mu
-\frac{\gamma}{2}f}\!\Big),\cn \quad \alpha_2:=\arctan\!\Big(\!
\frac\gamma 2\!\Big),\quad \alpha_3:=\eta^{-1}(\gamma)
\end{equation}
\cblue where
\begin{equation}
\label{fdef}
f:=\sqrt{(\gamma^2+4)(\mu-\mu_c)^2-4\mu^2}.
\end{equation}
The {\it critical value of $\mu_c$} is set as
\begin{equation}
\label{muccrit}
\mu_c^{\rm crit}:=\mu\Big[1-\frac{2}{\sqrt{\gamma^2+4}}\Big].
\end{equation}
For later use we observe the relationship
\begin{equation}
\label{a12}
\alpha_1^\pm=\arctan\!\Big(\!\frac{\gamma}{2}\!\Big)\pm
\arctan\!\Big(\frac{f}{2\mu}\Big)=\alpha_2\pm
\arctan\!\Big(\frac{f}{2\mu}\Big).
\end{equation}

\begin{Remark}
\label{rem3}
(i) Let $4\mu\not=\gamma f$ and $\mu_c\le\mu_c^{\rm crit}$ for $\mu>\mu_c$ or
$\mu_c\ge\mu\big(1+\frac{2}{\sqrt{\gamma^2+4}}\big)$. Then
$\alpha_1^-$, $\alpha_1^+$ are well-defined and solve
\begin{equation}
\label{f2}
\gamma\sa+2\ca=\frac{2\mu}{\mu-\mu_c}.
\end{equation}
(ii) For $\gamma<2$, the constant $\alpha_3$ has the alternative representation
\begin{equation}
\label{a3def2}
\alpha_3=\arctan\Big(\frac{4\gamma}{4-\gamma^2}\Big).
\end{equation}
\noindent{\bf Proof.} (i) With $S:=\sin(\alpha)$, $C:=\cos(\alpha)$,
$g:=\frac{2\mu}{\mu-\mu_c}$, we find $2C=g-\gamma S$ implying
$4C^2=4-4S^2=(g-\gamma S)^2$. This leads to the quadratic equation
\[ S^2-\frac{2\gamma g}{\gamma^2+4}S+\frac{g^2-4}{\gamma^2+4}=0 \]
with the solutions
\begin{equation}
\label{Spm}
S_\pm \;=\; \frac{\gamma g\pm\sqrt{\gamma^2g^2+(4-g^2)(\gamma^2+4)}}
{\gamma^2+4}=\frac{\gamma g\pm 2d}{\gamma^2+4}
\end{equation}
where $d:=(\gamma^2+4-g^2)^{1/2}$.
Eqn.~(\ref{Spm}) and $2C=g-\gamma S$ imply
\[ C_\pm=\frac{2g\mp\gamma d}{\gamma^2+4} \]
and (\ref{f2}) has the solutions
\begin{equation}
\label{sol1}
\alpha_1^-=\arctan\Big(\frac{S_-}{C_-}\Big),\qquad
\alpha_1^+=\arctan\Big(\frac{S_+}{C_+}\Big).
\end{equation}
Writing
\[ d=\frac{\big((\gamma^2+4)(\mu-\mu_c)^2-4\mu^2\big)^{1/2}}{(\mu-\mu_c)}
=:\frac{f}{\mu-\mu_c} \]
we arrive at the defining identities (\ref{mindef}) of $\alpha_1^\pm$.
The existence of $\alpha_1^\pm$ in (\ref{sol1}) requires for $\mu>\mu_c$
\begin{equation}
\label{ineq}
\gamma^2+4\ge g^2=\frac{4\mu^2}{(\mu-\mu_c)^2}\;\Longleftrightarrow\;
\mu_c^{\rm crit}:=\mu\Big(1-\frac{2}{\sqrt{\gamma^2+4}}\Big)\ge\mu_c.
\end{equation}
(ii) In view of (\ref{s0}), the equation $\eta(\alpha_3)=\gamma$ is equivalent
to
\begin{equation}
\label{a3eq}
\gamma\sin(\alpha_3)+2\cos(\alpha_3)=2.
\end{equation}
The calculus in (i) with $g:=2$ entails $d=\gamma$ and $\alpha=0$,
$\alpha=\arctan\Big(\frac{2\gamma\pm2\gamma}{4-\gamma^2}\Big)=\alpha_3$
as solutions of Eqn.~(\ref{a3eq}).
\end{Remark}

As minimizer of the Cosserat problem, $\alpha\equiv0$ is only present if
$\mu_c=0$, but then $0=\alpha_1^-$ which is why the zero solution is not part
of (\ref{mindef}). \cn

\begin{Corollary}
\label{cor1}
Any solution $u\in X:=\{W^{2,2}(\Omega)\;|\;u(0)=0,\,u(1)=\gamma\}$ to
(\ref{Wmin}) is monotonically increasing. The homogeneous function
$\uhom(x):=\gamma\, x$ solves the Euler-Lagrange equations
(\ref{S1}), (\ref{S2}).
Depending on the values of $\mu$ and $\mu_c$, the corresponding solution
$\alpha$ to (\ref{S1}), (\ref{S2}) is given by
\begin{align}
\label{ELsol1}
\mbox{(i) }\mu=\mu_c:\quad &
\alpha(x)\equiv\alpha_2.\\
& (\uhom,\alpha_2)\mbox{ is a local minimizer of }E.\nn\\
\label{ELsol2}
\mbox{(ii) }\mu_c=0:\quad &
\alpha(x)\equiv\cblue\alpha_1^-=0\cn,\;
\alpha(x)\equiv\alpha_2,\;\alpha(x)\equiv\alpha_3.\\
& (\uhom,\cblue0\cn)\mbox{ and }(\uhom,\alpha_3)
\mbox{ are local minimizers of }E.\nn\\
\label{ELsol3}
\mbox{(iii) } \mu_c>0,\,\mu\not=\mu_c:\quad &
\alpha(x)\equiv\alpha_2
\cblue,\; \alpha(x)\equiv\alpha_1^\pm \mbox{ whenever }\alpha_1^\pm
\mbox{ exist.}\cn\\
& \mbox{For } \cblue\mu_c>\mu_c^{\rm crit}\cn, (\uhom,\alpha_2)
\mbox{ is a local minimizer of }E.\nn\\
& \cblue(\uhom,\alpha_1^\pm)\mbox{ are local minimizers of }E\mbox{ if }
\mu>\mu_c,\; \mu_c\le\mu_c^{\rm crit}\mbox{ and }\alpha_1^\pm\not=\alpha_2.
\cn\nn
\end{align}
\end{Corollary}

\noindent{\bf Proof.} Due to Lemma~\ref{lem3}, any minimizer
$u\in W^{2,2}(\Omega)$ to (\ref{Wmin}) must be piecewise linear.
Choosing a function $u$ with $u(0)=0$, $u(1)=\gamma$ which is not
monotonically increasing enlarges the component $\frac{\mu}{2}|u'|^2$ in $W$,
cf. Eqn.~(\ref{Edef2}).
With (\ref{uzero}), this demonstrates the optimality of $\uhom(x)$ in
the class $X$ of regular solutions.
It remains to find the optimal values of $\alpha$, and the strict
positivity of the second variation $D^2_\alpha(\uhom,\alpha)$ is sufficient
for that, see (\ref{posdef}) below.

\vspace*{2mm}
For fixed $u\in X$ and a test function $\delta\alpha\in C^\infty((0,1);\,\R)$,
the second variation of $E$ with respect to $\alpha$ is
\begin{equation}
\label{2var1}
D_\alpha^2 E(u,\alpha)(\delta\alpha,\delta\alpha)
=\int_0^1\!\!\Big(W_{\alpha\alpha}(u',\alpha,\alpha')
-\ddx W_{\alpha\alpha'}(u',\alpha,\alpha')\Big)\delta\alpha(x)^2
+W_{\alpha'\alpha'}(u',\alpha,\alpha')\delta\alpha'(x)^2\dx
\end{equation}
where subscripts denote partial derivatives. Here we have,
cf. Eqn.~(\ref{Edef2}),
\[ W(u',\alpha,\alpha')=\frac{\mu}{2}|u'|^2+\frac{\mu}{2}\Big(
\sa u'-4\shalf\Big)^2+\frac{\mu_c}{2}\Big(\ca u'-2\sa\Big)^2 \]
such that $W_{\alpha\alpha'}(u',\alpha,\alpha')=W_{\alpha'\alpha'}
(u',\alpha,\alpha')=0$. So the second variation (\ref{2var1}) w.r.t. $\alpha$
simplifies to
\begin{equation}
\label{2var}
D_\alpha^2 E(u,\alpha)(\delta\alpha,\delta\alpha)=\int_0^1
W_{\alpha\alpha}(u',\alpha)\,\delta\alpha(x)^2\dx.
\end{equation}
Here and below we simply write $W(u',\alpha)$ instead of $W(u',\alpha,\alpha')$
due to $L_c=0$.

Direct computations reveal
\begin{align}
\label{faa}
\begin{split}
W_{\alpha}(u',\alpha) \;=\;\; & \mu\Big(\sa u'\!-\!4\shalf\Big)
\Big(\ca u'\!-\!2\sa\Big)\\
& \quad -\mu_c\Big(\ca u'\!-\!2\sa\Big)\Big(\sa u'\!+\!2\ca\Big),\\
W_{\alpha\alpha}(u',\alpha) \;=\;\; & (\mu\!-\!\mu_c)\Big(
\ca u'\!-\!2\sa\Big)^2+\mu_c\Big(\sa u'\!+\!2\ca\Big)^2\\
& \quad -\mu\Big(\sa u'\!-\!4\shalf\Big)\Big(\sa u'\!+\!2\ca\Big).
\end{split}
\end{align}

(i) For $\mu=\mu_c$, by direct investigation of (\ref{EL1a}), (\ref{EL1b})
we can verify that $(\uhom,\alpha_2)$ solves the Euler-Lagrange equations.
With (\ref{faa}) we find
\begin{align}
W_{\alpha\alpha}(\uhom',\alpha_2)=\; & \mu\Big(\sin(\alpha_2)
\gamma\!+\!2\cos(\alpha_2)\Big)^2-\mu\Big(\!\sin(\alpha_2)
\gamma\!-\!4\sin^2\!\big(\frac{\alpha_2}{2}\big)\Big)\Big(\!\sin(\alpha_2)
\gamma\!+\!2\cos(\alpha_2)\!\Big)\nn\\
\label{waa}
=\; & \mu\Big(\sin(\alpha_2)\gamma\!+\!2\cos(\alpha_2)\Big)\!\Big(\sin(\alpha_2)
\gamma+2\cos(\alpha_2)-\sin(\alpha_2)\gamma+\underbrace{4\sin^2\!\big(
\frac{\alpha_2}{2}\big)}_{\;=\;2-2\cos(\alpha_2)}\Big)\\
=\; & 2\mu\Big(\sin(\alpha_2)\gamma+2\cos(\alpha_2)\Big)\nn\\
=\; & (\gamma^2+4)\mu\cos(\alpha_2).\nn
\end{align}
Next we observe the identities
\begin{align}
\label{costan}
\cos(\arctan(t)) \;=\; &\frac{1}{\sqrt{t^2+1}},\qquad t\in\R,\\
\label{sintan}
\sin(\arctan(t)) \;=\; &\frac{t}{\sqrt{t^2+1}},\qquad t\in\R.
\end{align}
Consequently
\begin{align}
\label{cosa2}
\cos(\alpha_2) \;=\; & \cos\!\Big(\!\arctan\!\Big(\frac{\gamma}{2}\Big)\Big)=
\frac{2}{(\gamma^2\!+\!4)^{1/2}},\\
\label{sina2}
\sin(\alpha_2) \;=\; & \sin\!\Big(\!\arctan\!\Big(\frac{\gamma}{2}\Big)
\Big)=\frac{\gamma}{(\gamma^2\!+\!4)^{1/2}}
\end{align}
such that
\[ W_{\alpha\alpha}(\uhom',\alpha_2)=(\gamma^2\!+\!4)\mu\cos(\alpha_2)
=2\mu\sqrt{\gamma^2\!+\!4}. \]
With (\ref{2var}), this implies the strict positivity of the second variation
w.r.t. $\alpha$, i.e. there is a constant $K>0$ such that
\begin{equation}
\label{posdef}
D_\alpha^2 E(\uhom,\alpha_2)(\delta\alpha,\delta\alpha)\ge K\|\delta\alpha\|^2
\end{equation}
for any test function $\delta\alpha$, proving that $(\uhom,\alpha_2)$
is indeed a local minimizer of $E$, see, e.g. \cite{GH04}.

\vspace*{2mm}
(ii) For $\mu_c=0$, the Euler-Lagrange equations read
\begin{align}
\label{EL2a}
\Big[1-\frac{\cca}{2}\Big]u'' \;&=\; \big[\ca-\sa\ca u'+\ssa-\cca\big]\alpha',\\
\label{EL2b}
0 \;&=\; \big(\ca u'-2\sa\big)\Big(\sa u'-4\shalf\Big).
\end{align}
\cblue The definition (\ref{fdef}) implies $f=\gamma\mu$ for $\mu_c=0$, leading
to $\alpha_1^-=0$. \cn
Direct investigation of (\ref{EL2a}), (\ref{EL2b}) shows that
$(\uhom,\cblue0\cn)$, $(\uhom,\alpha_2)$ and $(\uhom,\alpha_3)$
are solutions of the Euler-Lagrange equations. With (\ref{faa}) we find 
\begin{align*}
W_{\alpha\alpha}(\uhom',\cblue0\cn) \;\;=\;\; & \mu\gamma^2>0,\\
W_{\alpha\alpha}(\uhom',\alpha_3) \;\;=\;\; & \mu\Big(\cos(\alpha_3)\gamma
-2\sin(\alpha_3)\Big)^2.
\end{align*}
Hence $W_{\alpha\alpha}(\uhom',\alpha_3)\!=\!0$ only if
$\cos(\alpha_3)\gamma\!=\!2\sin(\alpha_3)$ or equivalently
$\alpha_3\!:=\!\eta^{-1}(\gamma)=\arctan\!\big(\frac{\gamma}{2}\big)$, i.e.
only if $\gamma=\eta\big(\!\arctan(\gamma/2)\big)$, i.e. only if
\begin{equation}
\label{faapos}
\gamma=\frac{4\sin^2\big(\frac{\arctan(\gamma/2)}{2}\big)}
{\sin(\arctan(\gamma/2))}
=\frac{2-2\cos(\arctan(\gamma/2))}{\sin(\arctan(\gamma/2))}
\end{equation}
where (\ref{s0}) was used.
With (\ref{costan}) and (\ref{sintan}), Eqn.~(\ref{faapos}) becomes
\[ \gamma=\frac{2-\frac{2}{\sqrt{1+\gamma^2/4}}}{\frac{\gamma/2}
{\sqrt{1+\gamma^2/4}}}=\frac{2\sqrt{1+\frac{\gamma^2}{4}}-2}{\gamma/2}\;
\Longleftrightarrow\;\gamma^2=4\sqrt{1+\frac{\gamma^2}{4}}-4\;
\Longleftrightarrow\;\gamma^2+4=2\sqrt{4+\gamma^2}. \]
The last equality is only satisfied for $\gamma=0$, proving
$W_{\alpha\alpha}(\uhom',\alpha_3)>0$ for $\gamma\not=0$ which yields as in (i)
the strict positivity (\ref{posdef}) of the second variation w.r.t. $\alpha$
such that $(\uhom,\cblue\alpha_1^-\cn)$ and $(\uhom,\alpha_3)$
are local minimizers of $E$.

\vspace*{2mm}
(iii) For $\mu>0$, $\mu_c>0$, $\mu\not=\mu_c$, Lemma~\ref{lem3}, Case~1 in
(ii) shows that $(\uhom,\alpha_2)$ solves the Euler-Lagrange equations
(\ref{S1}), (\ref{S2}). \cblue As for the second factor in Eqn.~(\ref{S2}),
\[ (\mu-\mu_c)\big(\sa\gamma+2\ca-2\big)-2\mu_c=0\;\Longleftrightarrow\;
\sa\gamma+2\ca=\frac{2\mu}{\mu-\mu_c}. \]
Hence, by virtue of Remark~\ref{rem3}, $(\uhom,\alpha_1^-)$ and
$(\uhom,\alpha_1^+)$ solve the Euler-Lagrange equations (\ref{S1}), (\ref{S2})
whenever $\alpha_1^\pm$ exist. \cn 

\newpage
In order to \cblue investigate when \cn $(\uhom,\alpha_2)$ is a local
minimizer of $E$, we once more \cblue inspect \cn the second variation
(\ref{2var}) of $E$ w.r.t. $\alpha$. Starting from (\ref{faa}) we find
\begin{align}
W_{\alpha\alpha}(\uhom',\alpha_2) \;\;=&\;\; (\mu\!-\!\mu_c)\big(
\underbrace{\cos(\alpha_2)\gamma\!-\!2\sin(\alpha_2)}_{=\;0}\big)^2
+\mu_c\big(\sin(\alpha_2)\gamma\!+\!2\cos(\alpha_2)\big)^2\nn\\
& \; -\mu\Big(\sin(\alpha_2)\gamma
\underbrace{-4\sin^2\!\big(\frac{\alpha_2}{2}\big)}_{=\;2\cos(\alpha_2)-2}
\Big)\big(\sin(\alpha_2)\gamma+2\cos(\alpha_2)\big)\nn\\
=&\;\; (\mu_c\!-\!\mu)\big(\sin(\alpha_2)\gamma\!+\!2\cos(\alpha_2)\big)^2
+2\mu\big(
\underbrace{\sin(\alpha_2)\gamma+2\cos(\alpha_2)}_{=\;\frac{\gamma^2+4}{2}
\cos(\alpha_2)}\big)\nn\\
\label{fastep}
=&\;\; \frac{\gamma^2\!+\!4}{2}\cos(\alpha_2)\Big[2\mu+(\mu_c\!-\!\mu)\big(
\sin(\alpha_2)\gamma+2\cos(\alpha_2)\big)\Big].
\end{align}
Using the formulas (\ref{cosa2}), (\ref{sina2}) in (\ref{fastep}), this yields
\begin{align*}
W_{\alpha\alpha}(\uhom',\alpha_2) \;\;=&\;\;\frac{\gamma^2\!+\!4}{2}
\frac{2}{\sqrt{\gamma^2\!+\!4}}\Big[2\mu+(\mu_c\!-\!\mu)\Big(
\frac{\gamma^2}{\sqrt{\gamma^2\!+\!4}}+\frac{4}{\sqrt{\gamma^2\!+\!4}}\Big)
\Big]\\
=&\;\; \sqrt{\gamma^2\!+\!4}\Big[2\mu+(\mu_c\!-\!\mu)\sqrt{\gamma^2\!+\!4}
\Big].
\end{align*}
This shows the strict positivity (\ref{posdef}) of
$D_\alpha^2 E(\uhom,\alpha_2)(\delta\alpha,\delta\alpha)$ provided
\cblue $\mu_c>\mu_c^{\rm crit}$.\cn


\cblue
From Eqns.~(\ref{costan}), (\ref{sintan}),
\begin{align}
\label{scval1}
\gamma\sin(\alpha_1^-)+2\cos(\alpha_1^-) &= \frac{2\mu}{\mu-\mu_c},\quad
\gamma\cos(\alpha_1^-)-2\sin(\alpha_1^-)=\frac{f}{\mu-\mu_c},\\
\label{scval2}
\gamma\sin(\alpha_1^+)+2\cos(\alpha_1^+) &= \frac{2\mu}{\mu-\mu_c},\quad
\gamma\cos(\alpha_1^+)-2\sin(\alpha_1^+)=\frac{-f}{\mu-\mu_c}.
\end{align}
Owing to the identity (\ref{faa}), this yields
\begin{align}
W_{\alpha\alpha}(\gamma,\alpha_1^-)=W_{\alpha\alpha}(\gamma,\alpha_1^+)\; &=\;
(\mu-\mu_c)\frac{f^2}{(\mu-\mu_c)^2}+\mu_c\frac{4\mu^2}{(\mu-\mu_c)^2}
-\mu\Big(\frac{2\mu}{\mu-\mu_c}-2\Big)\frac{2\mu}{\mu-\mu_c}\nn\\
\label{Wconv}
&= \; \frac{f^2}{\mu-\mu_c}+\frac{4\mu^2\mu_c}{(\mu-\mu_c)^2}
-\frac{2\mu^2}{\mu-\mu_c}\frac{2\mu_c}{\mu-\mu_c}=\frac{f^2}{\mu-\mu_c}.
\end{align}
For $\mu>\mu_c$ and $f\not=0$, this proves that $(\uhom,\alpha_1^\pm)$ are
local minimizers of $E$ while for $f=0$, it holds
$\alpha_1^\pm=\arctan(\frac{\gamma}{2})=\alpha_2$
due to (\ref{a12}). For $\mu_c<\mu_c^{\rm crit}$, Eqn.~(\ref{ineq}) holds with
strict inequality such that
\begin{align}
W_{\alpha\alpha}(\gamma,\alpha_2) \;&=\; 2\mu(\gamma^2+4)^{1/2}+(\mu-\mu_c)\big[
0-(\gamma^2+4)\big]\nn\\
&=\; (\gamma^2+4)^{1/2}\big[2\mu-(\mu-\mu_c)(\gamma^2+4)^{1/2}\big]\nn\\
\label{Wnconv}
&<\; (\gamma^2+4)^{1/2}\big[2\mu-2\mu]=0
\end{align}
and $(\uhom,\alpha_2)$ is no minimizer of $E$ for $\mu_c<\mu_c^{\rm crit}$.
\cn
\vspace*{-10pt}\qed

\newpage

\begin{Remark}
\label{rem4}
We point out that the solution $(\uhom,\alpha_4)$ of the Euler-Lagrange
equations found in Lemma~\ref{lem3}, (ii) is not an independent case. 
The solution $(\uhom,\alpha_4)$ requires $\big|\frac{\mu-\mu_c}{\mu}\big|\le1$,
cf. Eqn.~(\ref{rend}), and then, by the very definition of Case~2 in
Lemma~\ref{lem3},
\begin{align}
\sin(\alpha_4)\gamma \;\;=\;\; & 4\sin^2\!\big(\frac{\alpha_4}{2}\big)
+\frac{2\mu_c}{\mu\!-\!\mu_c} \;\;=\;\; 2-2\cos(\alpha_4)
+\frac{2\mu_c}{\mu\!-\!\mu_c}\nn\\
\label{stp1}
\;\;=\;\; & 2-\frac{2(\mu\!-\!\mu_c)}{\mu}+\frac{2\mu_c}{\mu\!-\!\mu_c}
\;\;=\;\; \frac{2\mu_c}{\mu}+\frac{2\mu_c}{\mu\!-\!\mu_c}
\;\;=\;\; \frac{2\mu_c(2\mu\!-\!\mu_c)}{\mu(\mu\!-\!\mu_c)}.
\end{align}
Due to the relationship
\begin{equation}
\label{sincos}
\sin(\arccos(t))=\sqrt{1-t^2} \qquad \mbox{for }t\!\in\![-1,+1],
\end{equation}
this yields
\begin{equation}
\label{stp0}
\sin(\alpha_4)=\sqrt{1-\frac{(\mu-\mu_c)^2}{\mu^2}}=\frac{1}{\mu}
\sqrt{\mu_c(2\mu-\mu_c)}.
\end{equation}
With (\ref{stp1}) we obtain
\begin{equation}
\label{gammacond}
\gamma \;\;=\;\; \frac{2\mu_c(2\mu-\mu_c)}{\mu(\mu-\mu_c)\sin(\alpha_4)}
\;\;=\;\; \frac{2\mu_c(2\mu-\mu_c)}{\mu(\mu-\mu_c)}\;
\frac{\mu}{\sqrt{\mu_c(2\mu-\mu_c)}}
\;\;=\;\; \frac{2\sqrt{\mu_c(2\mu-\mu_c)}}{\mu-\mu_c}.
\end{equation}
For $\gamma$ given by (\ref{gammacond}), we have
\begin{align*}
\gamma\cos(\alpha_4)-2\sin(\alpha_4) \;\;=&\;\; \gamma\;\frac{\mu-\mu_c}{\mu}
-\frac{2}{\mu}\sqrt{\mu_c(2\mu-\mu_c)}\\
\;\;=&\;\; \frac{2\sqrt{\mu_c(2\mu-\mu_c)}}{\mu-\mu_c}\;\frac{\mu-\mu_c}{\mu}
-\frac{2}{\mu}\sqrt{\mu_c(2\mu-\mu_c)} \;\;=\;\; 0,
\end{align*}
leading to $\gamma\cos(\alpha_4)=2\sin(\alpha_4)$
or equivalently $\frac{\gamma}{2}=\tan(\alpha_4)$. This shows
\begin{equation}
\label{eqalpha}
\alpha_4=\arctan\!\Big(\frac{\gamma}{2}\Big)=\alpha_2
\end{equation}
for $\gamma$ given by (\ref{gammacond}).

In order to have $\gamma>0$ in (\ref{gammacond}), it must hold $0<\mu_c<\mu$,
and the condition $\big|\frac{\mu-\mu_c}{\mu}\big|\le1$
is automatically satisfied.
\end{Remark}

\begin{Remark}
\label{rem5}
The following table lists the mechanical energies corresponding to the
local minimizers of Corollary~\ref{cor1}.
\begin{align}
\label{Een2i}
\mbox{(i) }\mu=\mu_c:\quad &
E(\uhom,\alpha_2)=\mu\Big[\gamma^2+4-2\sqrt{\gamma^2\!+\!4}\Big].\\
\label{Een13}
\mbox{(ii) }\mu_c=0:\quad &
E(\uhom,\cblue\alpha_1^-\cn)=E(\uhom,\alpha_3)=\frac{\mu}{2}\gamma^2.\\
\label{Een2iii}
\mbox{(iii) } \mu_c>0,\,\mu\not=\mu_c:\quad &
E(\uhom,\alpha_2)=\mu\Big[\gamma^2+4-2\sqrt{\gamma^2\!+\!4}\Big]
\mbox{ if\/ } \cblue\mu_c>\mu_c^{\rm crit}\cn.\\
\label{Een2iiib}
& \cblue E(\uhom,\alpha_1^\pm)=\frac{\mu+\mu_c}{2}\gamma^2
-\frac{2\mu_c^2}{\mu-\mu_c}
\mbox{ whenever }\alpha_1^\pm \mbox{ exist and }\alpha_1^\pm\not=\alpha_2.
\cn
\end{align}
\end{Remark}

\begin{figure}[h!bt]
\unitlength1cm
\begin{picture}(11.5,3.4)
\put(-0.1,-0.40){\psfig{figure=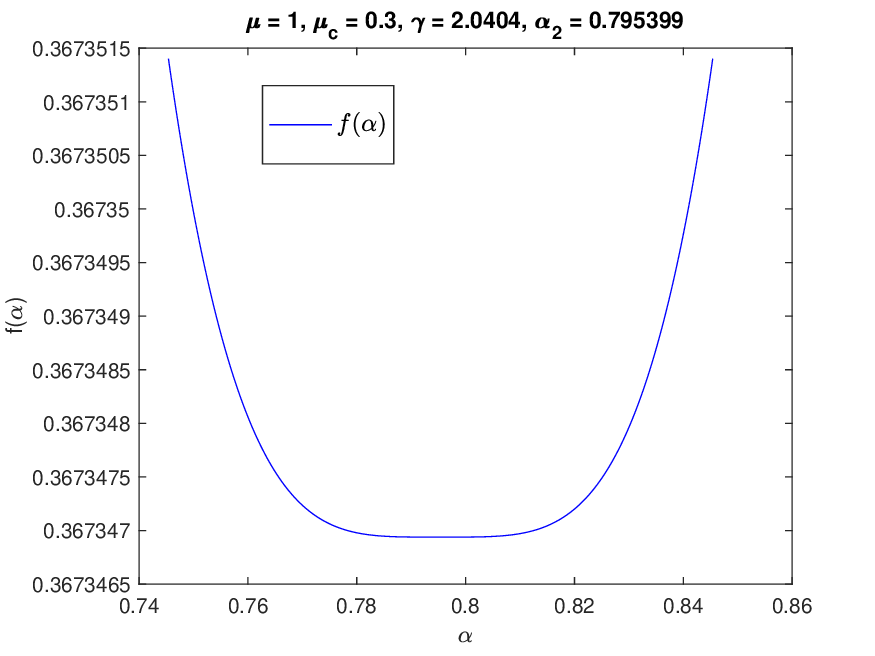,width=5.9cm}}
\put(4.8,-0.40){\psfig{figure=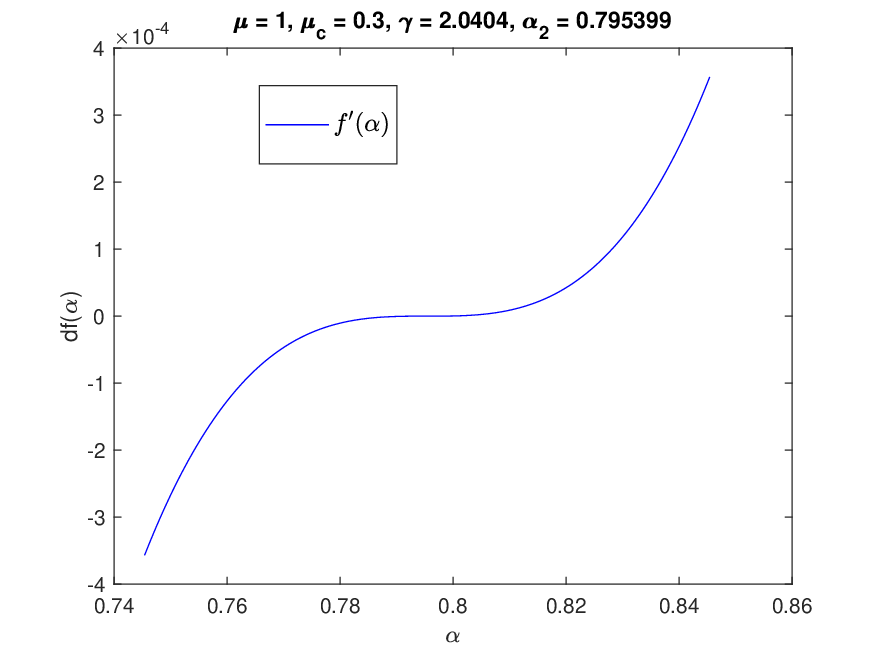,width=5.9cm}}
\put(10.5,-0.40){\psfig{figure=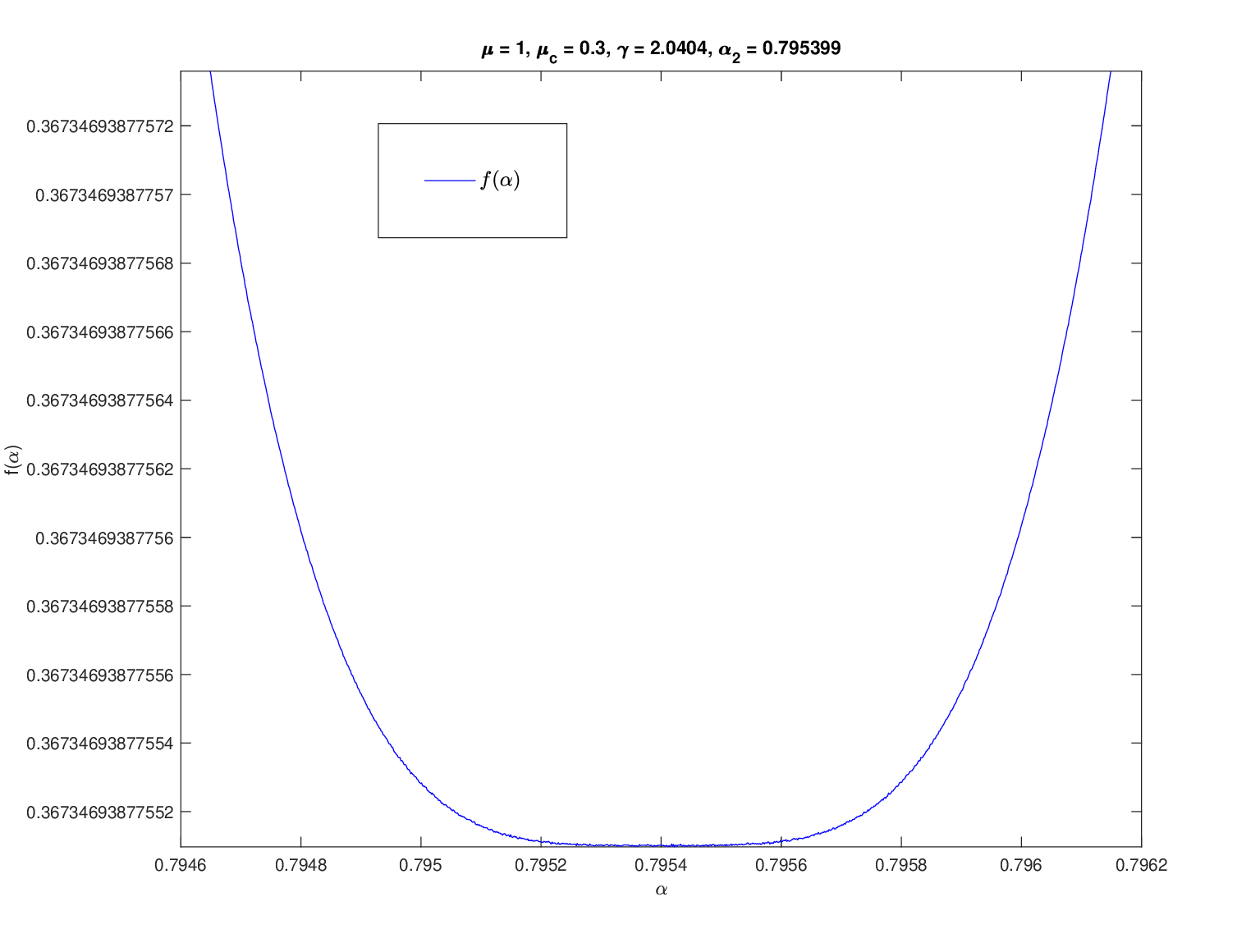,width=5.9cm}}
\end{picture}
\caption{\label{fig6a}
Left and center: Plot of $f(\alpha)$, $f'(\alpha)$ near $\alpha_2$ for
$\mu=1$, $\mu_c=0.3$ and $\gamma$ given by (\ref{gammacond}).
Right: Close-up of $f(\alpha)$ near $\alpha_2$. The tiny oscillations of the
graph near $\alpha_2$ in the close-up are numerical rendering artifacts
generated by MATLAB.}
\end{figure}

\vspace*{1mm}
\begin{Remark}
\label{rem6}
We want to investigate whether $(\uhom,\alpha_2)$ is a
minimizer of $E$ in the situation of Remark~\ref{rem4}, i.e. if
\begin{equation}
\label{gammadef}
0<\mu_c<\mu\mbox{ and }(\mu-\mu_c)\gamma=2\sqrt{\mu_c(2\mu-\mu_c)}.
\end{equation}
As $\uhom'(x)\equiv\gamma$ and assuming that $\alpha(x)\equiv a$ for
constant $a\in\R$ and $x\in(0,1)$,
it is enough to investigate the real function
\[ f(a):=W(\uhom,a)=\frac{\mu}{2}\gamma^2+\frac{\mu}{2}\Big(
\gamma\sin(a)-4\sin^2\!\big(\frac{a}{2}\big)\Big)^2+\frac{\mu_c}{2}
\Big(\gamma\cos(a)-2\sin(a)\Big)^2, \]
cf. the definition of $W$ in (\ref{Edef2}). Straightforward computations yield
\begin{align*}
f'(a) \;\;=&\;\; \big(\gamma\cos(a)-2\sin(a)\big)\big[(\mu-\mu_c)
\big(\gamma\sin(a)+2\cos(a)\big)-2\mu\big],\\
f''(a) \;\;=&\;\; (\mu-\mu_c)\big[\big(\gamma\cos(a)-2\sin(a)\big)^2
-\big(\gamma\sin(a)+2\cos(a)\big)^2\big]+2\mu\big(\gamma\sin(a)+2\cos(a)\big),\\
f'''(a) \;\;=&\;\; \big(\gamma\cos(a)-2\sin(a)\big)\big[2\mu-4(\mu-\mu_c)\big(
\gamma\sin(a)+2\cos(a)\big)\big],\\
f^{(4)}(a) \;\;=&\;\; 4(\mu-\mu_c)\big[\big(\gamma\sin(a)+2\cos(a)\big)^2-\big(
\gamma\cos(a)-2\sin(a)\big)^2\big]-2\mu\big(\gamma\sin(a)+2\cos(a)\big).
\end{align*}
Remarkably, $f'(\alpha_2)=f''(\alpha_2)=f'''(\alpha_2)=0$, but
\[ f^{(4)}(\alpha_2)=\frac{12\mu^2}{\mu-\mu_c}>0, \]
showing that $\alpha_2$ is a minimizer of $f$ and indicating that
$(\uhom,\alpha_2)$ is a local minimizer of $E$ for constant
$\alpha(x)\equiv\alpha_2$ under the choice of parameters
(\ref{gammadef}). The minimal energy coincides with (\ref{Een2iii}) for
$\gamma$ given by (\ref{gammacond}).
Due to $f'(\alpha_2)=f''(\alpha_2)=f'''(\alpha_2)=0$, the function $f$ is
extremely flat near $\alpha_2$, making it very difficult to numerically compute
the correct minimizer, see Fig.~\ref{fig6a}.
\end{Remark}

\vspace*{4mm}

\begin{Remark}
\label{rem7}
Let $\mu>0$ and $\gamma>0$ be given.
\cblue
Corollary~\ref{cor1} states that the value of $\mu_c$ selects the minimizer
of $E(\uhom,\cdot)$. For $\mu_c>\mu_c^{\rm crit}$, $(\uhom,\alpha_2)$ is
optimal while for $\mu_c<\mu_c^{\rm crit}$ and
$\mu>\mu_c$, $(\uhom,\alpha_1^{\pm})$ are local minimizers. These findings
are confirmed by our numerical simulations, cf. Fig.~\ref{fig10}.\cn
\end{Remark}

\begin{Remark}
\label{rem8}
It is instructive to compare Corollary~\ref{cor1} and Remark~\ref{rem5} with
the results in \cite{FN17}, stating that
\begin{align}
\label{scale1}
\mbox{For }\mu_c\ge\mu>0:&\qquad W_{\mu,\mu_c}(R;F)\sim W_{1,1}(R;F),\\
\label{scale2}
\mbox{For }\mu>\mu_c\ge0:&\qquad W_{\mu,\mu_c}(R;F)\sim W_{1,0}(R;
\widetilde{F}_{\mu,\mu_c}).
\end{align}
Here, $F=D\phi\in\R^{n\times n}$ is a given deformation gradient,
$\widetilde{F}_{\mu,\mu_c}:=\frac{\mu-\mu_c}{\mu}F$,
\[ f\sim g :\Longleftrightarrow \argmin_{R\in\SO(n)}f(R)=\argmin_{R\in\SO(n)}
g(R) \]
for two functions $f$ and $g$, and
\begin{equation}
\label{Wnew}
W_{\mu,\mu_c}(R;F):=\mu\big|\Sym(R^TF-\id_3)\big|^2
+\mu_c\big|\skw(R^TF-\id_3)\big|^2
\end{equation}
which coincides with (\ref{Edef2}) for $L_c=0$,
$F$ given by (\ref{Fdef}), and $R\in\SO(3)$ defined by (\ref{R}).

\cblue In the classical parameter range $\mu_c\ge\mu>0$,
\begin{align*}
W_{1,1}(R;F) \;&=\; \frac12\Big[\gamma^2+(\gamma\sa+2\ca-2)^2+(\gamma\ca-2\sa)^2
\Big]\\
&\sim\;\Big(\gamma\sa+2\ca-2\big)\Big)^2+\Big(\gamma\ca-2\sa\Big)^2.
\end{align*}
The value $\alpha_3=\eta^{-1}(\gamma)$ corresponds to the zeros of the first
quadratic term on the right, the global minimizer $\alpha_2$ to the zeros of
the second quadratic term. This includes the case $\mu=\mu_c$.

Similarly, for the non-classical parameter range $\mu>\mu_c\ge0$, introducing
the scaling parameter $r=r_{\mu,\mu_c}:=\frac{\mu-\mu_c}{\mu}$,
direct inspection reveals
\begin{align}
W_{1,0}(R;\widetilde{F}) \;=&\;
\frac12\Big(\gamma r\sa+2r\ca-2\Big)^2+\frac12\gamma^2 r^2+(r-1)^2\nn\\
\label{Wsim}
\;\sim&\; \Big(\gamma\sa+2\ca-\frac{2}{r}\Big)^2=\Big(\gamma\sa+2\ca
-\frac{2\mu}{\mu-\mu_c}\Big)^2\mbox{ for }\mu\not=\mu_c.
\end{align}
The zeros on the right are the global minimizers $\alpha_1^\pm$ whenever they
exist.
Especially, for $\mu_c=0$, Eqn.~(\ref{Wsim}) with the help of (\ref{s0}) states
that
\[ W_{1,0}(R;\widetilde{F})\sim
\Big(\gamma\sa-4\sin^2\!\big(\frac{\alpha}{2}\big)\Big)^2. \]
As outlined in Remark~\ref{rem3}, the right hand side is zero (and minimal) for
$\alpha=0=\alpha_1^-$ and for $\alpha=\alpha_3$
corresponding to (\ref{ELsol2}).

\vspace*{2mm}
\cn In conclusion, Corollary~\ref{cor1} and Remark~\ref{rem5} confirm
the results in \cite{FN17}.

\end{Remark}

\section{Numerical simulations for vanishing internal length scale}
\label{secnum}
For $L_c>0$, various advanced numerical tools such as multigrid methods
are available, \cite{Ble20}. We will not discuss this here. In contrast,
the case $L_c=0$ is numerically challenging as the regularizing term
$2\mu L_c^2|\alpha'|^2$ in $E$ and $\Ered$ disappears in the limit.
We investigate the problem in two different ways and compare the solution
strategies. We begin with the situation of non-regular solutions with
$\alpha\in L^4(\Omega)$.

\vspace*{4mm}
\begin{enumerate}
\item[\textcircled{1}]\underline{Newton-GMRES algorithm}

For $L_c=\mu_c=0$ and $\alpha\in L^4(\Omega)$, the solutions $(u,\alpha)$ to
$\Ered(u,\alpha)\to\min$ satisfy pointwise for $x\in\Omega=(0,1)$ the two
purely algebraic equations (see \cite[Eqn.~(3.30)]{neff2009simple}
for a derivation),
\begin{align}
\label{Lim1}
\alpha(\alpha-u')\Big(\alpha-\frac{u'}{2}\Big) &=0,\\
\label{Lim2}
(u'-\zeta)\Big(u'+\frac18(u')^3-\zeta\Big) &=0
\end{align}
subject to the boundary conditions
\begin{equation}
\label{Lim3}
u(0)=0,\quad u(1)=\gamma,\quad u'(0)=u'(1),\quad \alpha(0)=\alpha(1)=\alpha_D
\end{equation}
for a constant $\alpha_D\in\R$ and with $\zeta:=(1+\alpha_D^2)u'(0)-\alpha_D^3$.

\vspace*{2mm}
In general, there is a multitude of solutions to (\ref{Lim1})--(\ref{Lim3}).
In order to specify a solution, one may prescribe a volume fraction
$\theta\in[0,1]$ such that
\begin{equation}
\label{mass}
\cL^1\big(\{x\in\Omega\;|\;\alpha(x)=0\}\big)=\theta,
\end{equation}
where $\cL^1$ denotes the one-dimensional Lebesgue measure.

\vspace*{2mm}
In order to compute solutions to (\ref{Lim1})--(\ref{Lim3}),
the derivative $u'$ is approximated by central difference quotients
leading to a discrete problem in the standard form
\begin{equation}
\label{zero}
\mbox{Find }x\in\R^{2N-3}: \quad G(x)=0
\end{equation}
for a differentiable function $G:\R^{2N-3}\mapsto\R^{2N-3}$ and with
$N\in\N$ denoting the number of discretization points in $\Omega$.
The solutions to (\ref{zero}) are computed using a Newton-GMRES method
where the Fr{\'e}chet derivative $DG$ of $G$ is approximated by
\begin{equation}
\label{DG}
DG(x)d\sim\frac{G(x+\delta d)-G(x)}{\delta}
\end{equation}
for suitable small $\delta>0$ and GMRES is used to solve the linearized
equations. While the algorithm computes $DG$ automatically thanks to
Eqn.~(\ref{DG}), its practicability is limited by the huge memory requirements.

The implementation details related to the definition of $G$ are left
out here.

\vspace*{3mm}
\item[\textcircled{2}]\underline{BFGS Quasi-Newton method}

As in the first method, the values of $\alpha$, $\alpha'$ and $u'$
are discretized along $N\in\N$ points $z_i\in\Omega$. Letting
$f_i:=W(u'(z_i),\alpha(z_i),\alpha'(z_i))$ be the integrand either in
(\ref{Wdef}) or in (\ref{Wreddef}), the repeated Trapezoidal rule
\[ \Bigg(\frac12 f_1+\sum_{i=2}^{N-1}f_i+\frac12 f_N\Bigg) \]
is used for the approximate integration. The minimization of the functional
is carried out by a Quasi-Newton method where the approximate
Hessian is computed by the Broyden-Fletcher-Goldfarb-Shanno update formula,
\cite{B70}. Again, we leave out the implementation details, but refer to
\cite{Ble13,Ble14,Ble17}, where the limited-memory variant of the algorithm is
applied to Cosserat plasticity and recrystallization.
\end{enumerate}

\vspace*{1mm}
In comparison, both algorithms compute the same solutions. However,
the BFGS Quasi-Newton method is capable of handling larger values of $N$
due to the tremendous memory needs of the Newton-GMRES scheme.

\begin{figure}[pht]
\unitlength1cm
\begin{picture}(11.5,17.0)
\put(3.5,-0.50){\psfig{figure=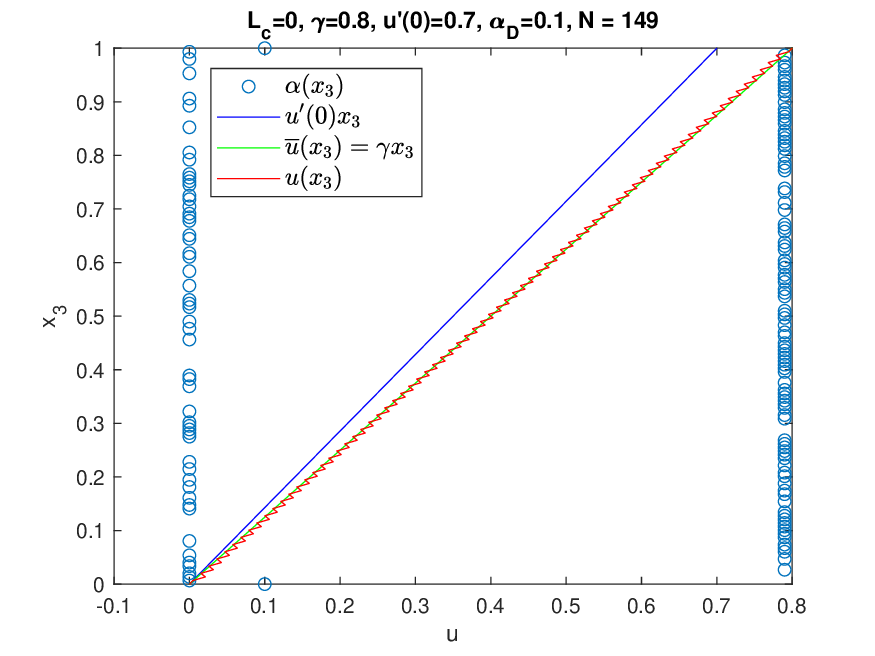,width=8.3cm}}
\put(3.5,7.00){\psfig{figure=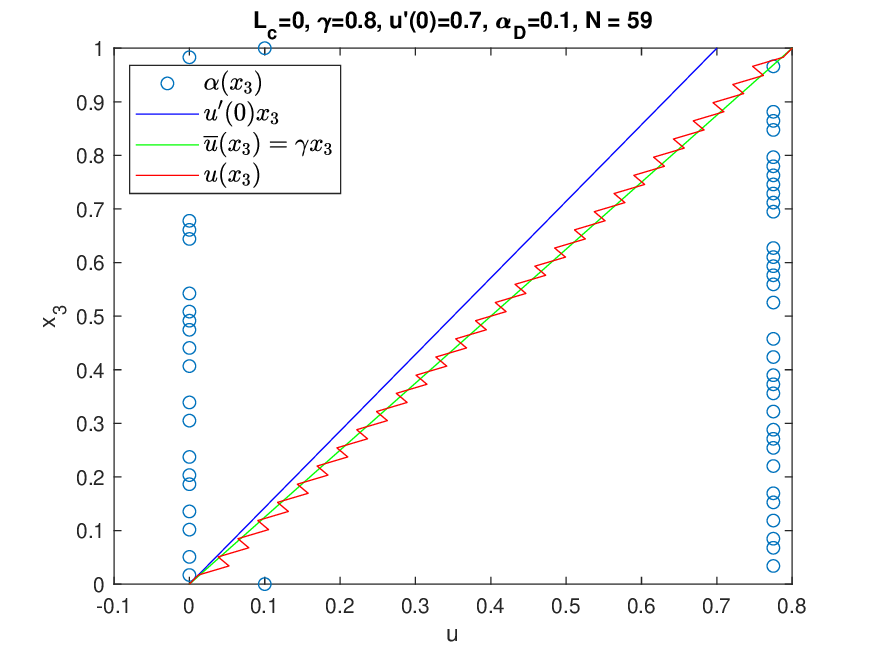,width=8.3cm}}
\put(3.5,14.00){\psfig{figure=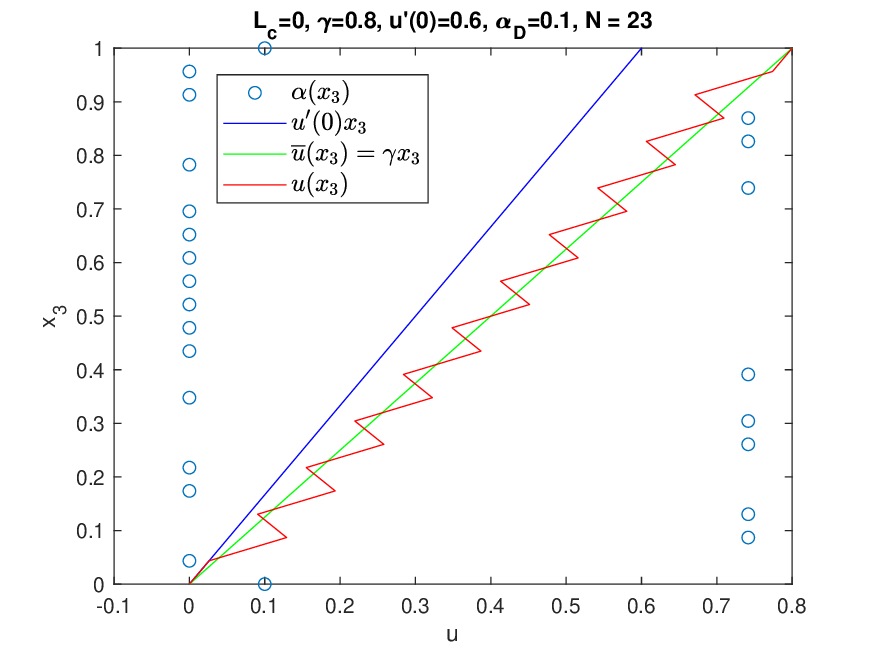,width=8.3cm}}
\end{picture}
\caption{\label{fig5}
Minimizers of $\Ered$ in $\cX^D$ for $N=23$, $N=59$ and $N=149$ 
for $\mu=1$, $\mu_c=0$, $\gamma=0.8$, prescribed $u'(0)=u'(1)$
(blue line) and $\alpha_D=u(0)=u(1)=0.1$.
The blue balls are the computed values of
$\alpha$ randomly concentrating at $0$ and $\gamma$.
The green line displays the homogeneous deformation
$\overline{u}(x)=\gamma\, x$. The deformation $u$ is
rendered in red forming a microstructure with a sawtooth pattern.}
\end{figure}

\vspace*{3mm}
Fig.~\ref{fig5} displays three minimizers of $\Ered$ in $\cX^D$ for different
values of $N$ and prescribed slope $u'(0)$ at the boundary.
The deformation $u_N$ forms a sawtooth pattern with alternating, constant
slopes, while the values of $\alpha$ randomly concentrate at $0$ and $\gamma$.
The computations of Fig.~\ref{fig5} suggest further that $u_N$ converges to
the homogeneous deformation $\overline{u}(x):=\gamma\,x$ as $N\to\infty$.

\begin{figure}[pht]
\unitlength1cm
\begin{picture}(11.5,15.9)
\put(3.5,-0.10){\psfig{figure=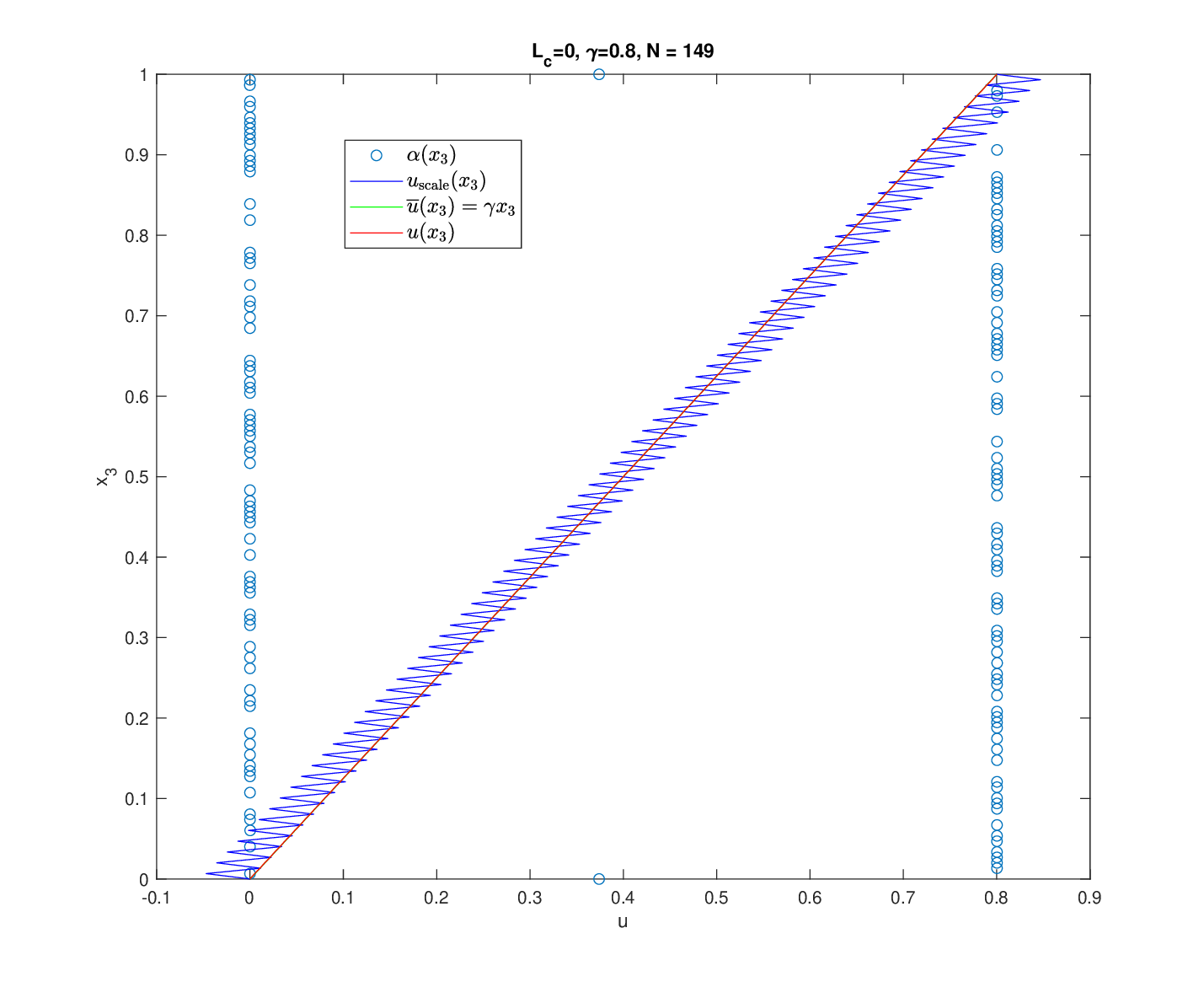,width=8.3cm}}
\put(3.5,6.80){\psfig{figure=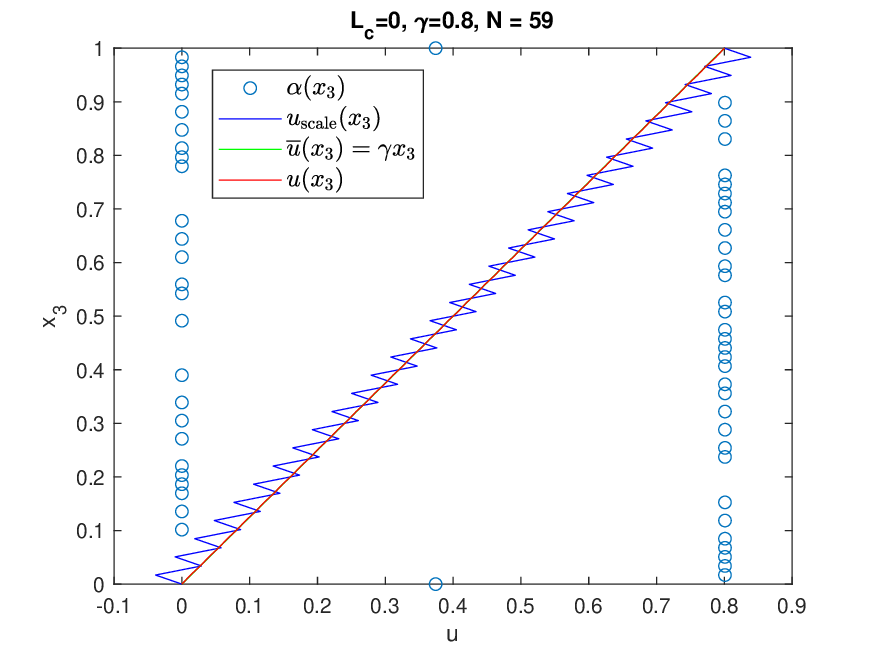,width=8.3cm}}
\put(3.5,13.20){\psfig{figure=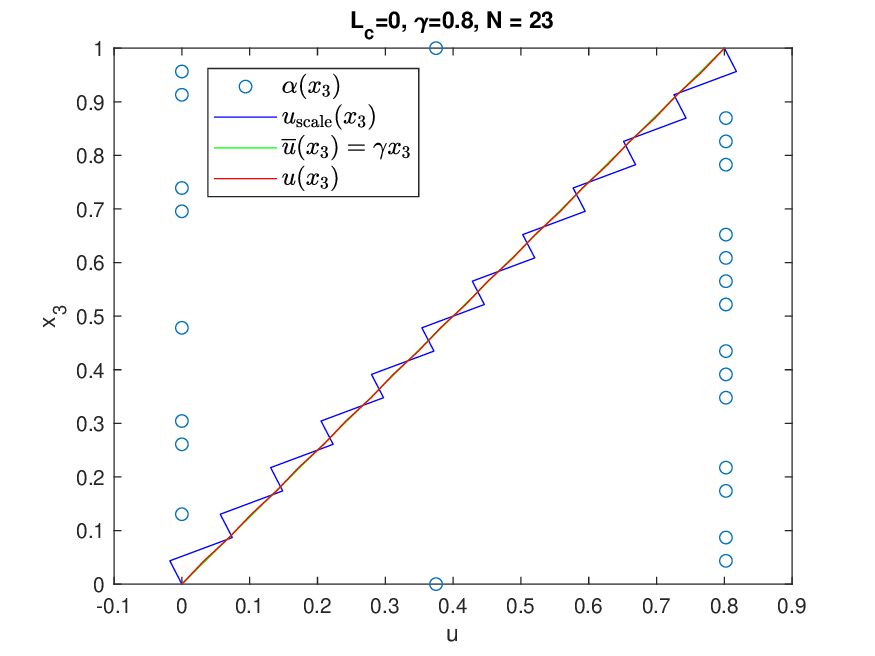,width=8.3cm}}
\end{picture}
\caption{\label{fig6}
Minimizers of $\Ered$ in $\cXCred$ for $N=23$, $N=59$ and $N=149$
for $\mu=1$, $\mu_c=0$, $\gamma=0.8$.
The blue balls are the computed values of
$\alpha$ randomly concentrating at $0$ and $\gamma$.
The deformation $u$ is rendered in red and forms a microstructure with a
sawtooth pattern which is extremely close to the homogeneous deformation
$\overline{u}(x)=\gamma\, x$ rendered in green.
The differences between $u$ and $\overline{u}$ is amplified by the blue line
rendering the rescaled deformation
$u_\mathrm{scale}(x):=\gamma\, x+N(u(x)\!-\!\gamma x)$.
The value $\alpha(0)$ is strictly smaller than $\frac\gamma2$ and indicates
by which amount the sawtooth
solution $u$ energetically beats the homo\-geneous solution.}
\end{figure}

\begin{figure}[pht]
\unitlength1cm
\begin{picture}(11.5,16.3)
\put(3.5,+0.50){\psfig{figure=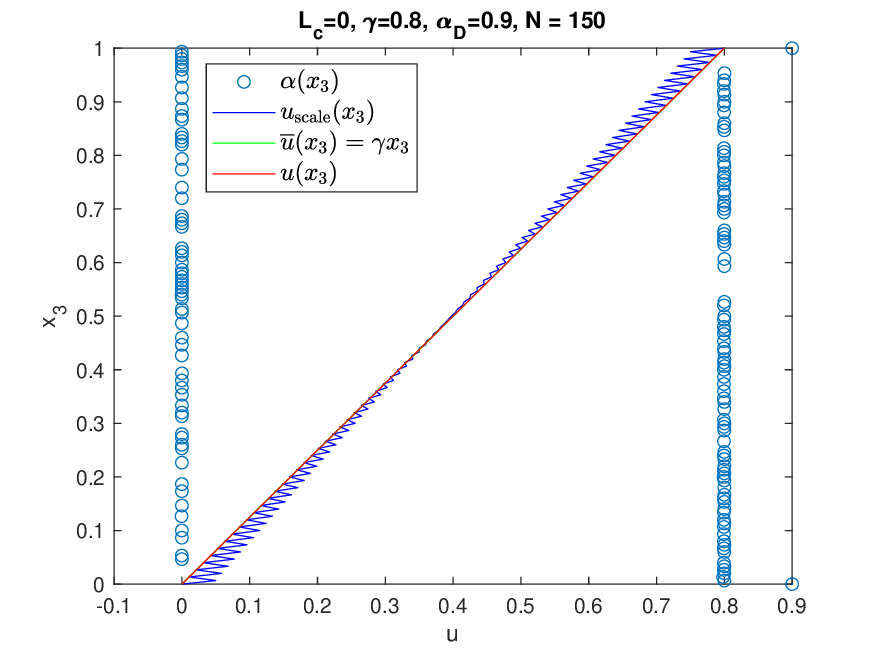,width=8.3cm}}
\put(3.5,7.40){\psfig{figure=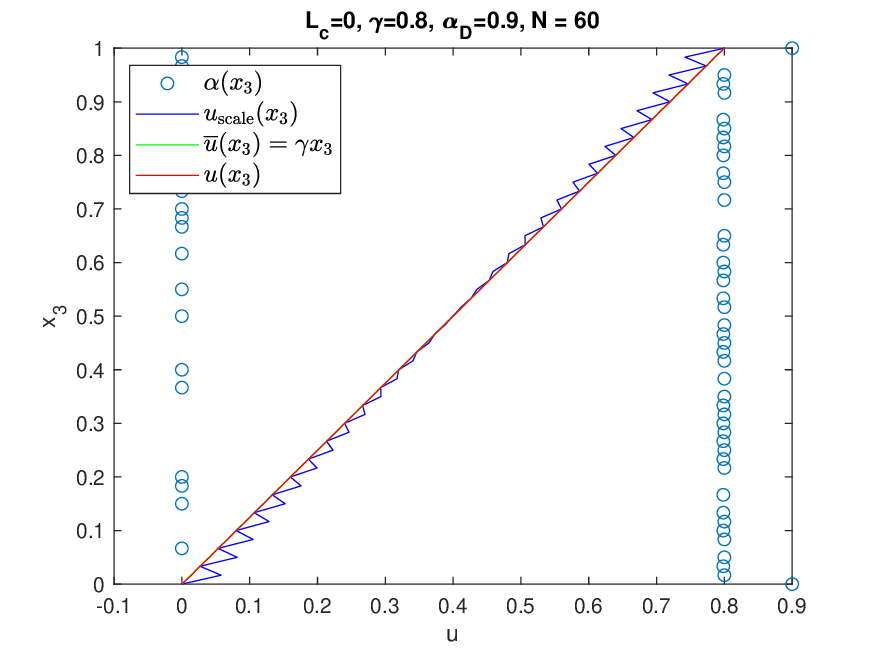,width=8.3cm}}
\put(3.5,13.80){\psfig{figure=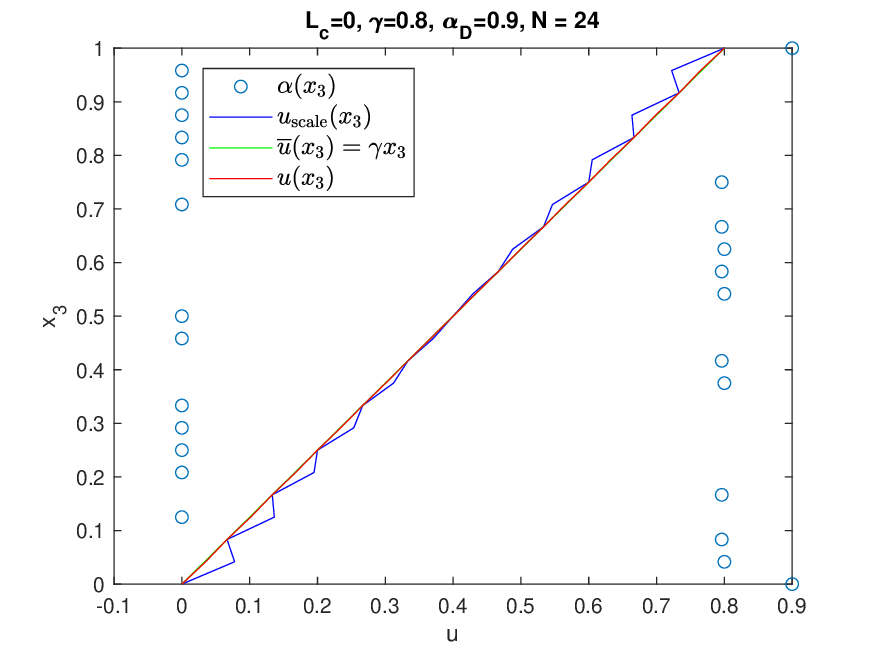,width=8.3cm}}
\end{picture}
\caption{\label{fig7}
Minimizers of $\Ered$ in $\cX^D$ for $N=24$, $N=60$ and $N=150$
for $\mu=1$, $\mu_c=0$, $\gamma=0.8$. The blue balls are the
computed values of $\alpha$ randomly concentrating at $0$ and $\gamma$.
The deformation $u$ is rendered in red and forms a microstructure with a
sawtooth pattern which is extremely close to the homogeneous deformation
$\overline{u}(x)=\gamma\, x$ rendered in green. The blue line is the rescaled
deformation $u_\mathrm{scale}(x):=\gamma\, x+N(u(x)\!-\!\gamma\, x)$ amplifying
the differences between $u$ and $\overline{u}$. As can be seen, the oscillation
of $u$ is largest near $\partial\Omega$ and smallest at $x=0.5$.}
\end{figure}

Fig.~\ref{fig6} and Fig.~\ref{fig7} compare the minimizers of $\Ered$ in
$\cX^C$ and $\cX^D$ for three different values of $N$.
Again, $\alpha$ randomly concentrates at $0$ and $\gamma$ while the deformation
$u_N$ forms a sawtooth pattern with alternating, constant slopes.
The computed deformation $u$ is extremely close to
$\overline{u}$ but energetically beats the homogeneous solution.

\vspace*{2mm}
We will adopt the following notation in the set $\cM^+(\R)$ of finite positive
{\em Radon measures}.
For a sequence $(v_j)_{j\in\N}$ and a Carath{\'e}odory function $f$ such that
\[ f(x,v_j)\rightharpoonup\Bigg(x\mapsto\int\limits_\R
f(x,A)\,\mathrm{d}\nu_x(A)\Bigg)\quad\mbox{in }L^1(\R) \]
for a parameterized measure $\nu=(\nu_x)_{x\in\R}$ we use the shorthand
notation $v_j\weakY\nu$.

\begin{Remark}
\label{rem9}
The numerical computations in Figs.~\ref{fig5}--\ref{fig7} indicate the
following result regarding non-regular solutions.

Let $L_c=0$ and $(u_N,\alpha_N)_{N\in\N}$ be a minimizing
sequence of $\Ered$ in $W^{1,2}(\Omega)\times L^4(\Omega)$ with
$u_N(x)\rightharpoonup\uhom(x)$, $\alpha_N\rightharpoonup\alpha$ for
$N\to\infty$. Then, there is a (not relabelled) subsequence and a constant
$\theta\in[0,1]$ such that
\begin{equation}
\label{Yconv}
\Wred(u_N,\alpha_N)\weakY\theta\delta_0+(1-\theta)\delta_\gamma.
\end{equation}
The existence of a limiting Young measure generated by the subsequence
$(u_N,\alpha_N)$ can be made rigorous by applying the fundamental theorem of
Young measures, see e.g. \cite{M98}. We observe that the family
$(\Wred(u_N,\alpha_N))_{N\in\N}$ is uniformly bounded in $L^1$ and
equi-integrable by the Dunford-Pettis theorem.
\end{Remark}

\begin{figure}[pht]
\unitlength1cm
\begin{picture}(11.5,16.3)
\put(1.5,11.00){\psfig{figure=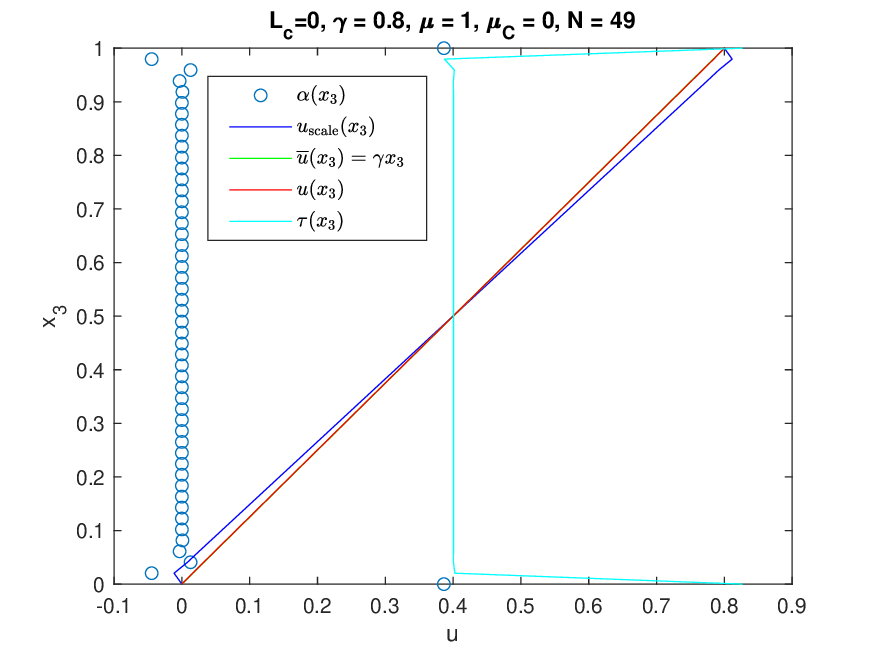,width=6.6cm}}
\put(7.5,11.00){\psfig{figure=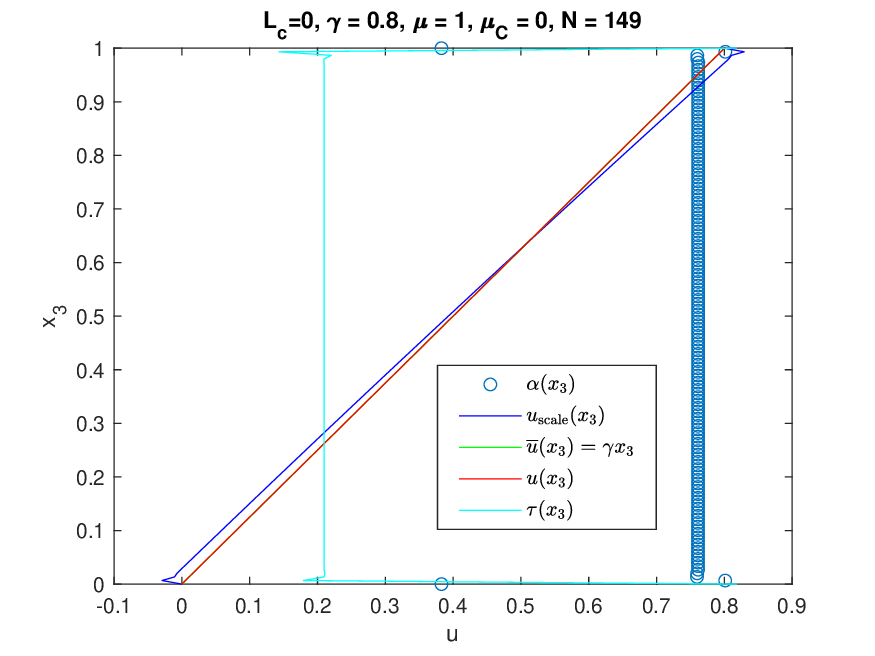,width=6.6cm}}
\put(1.5,5.50){\psfig{figure=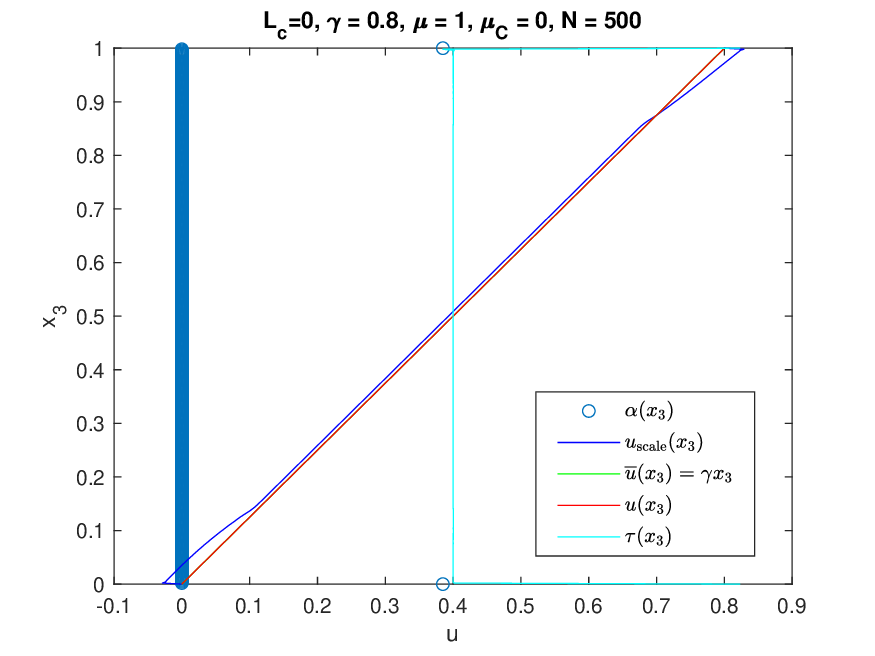,width=6.6cm}}
\put(7.5,5.50){\psfig{figure=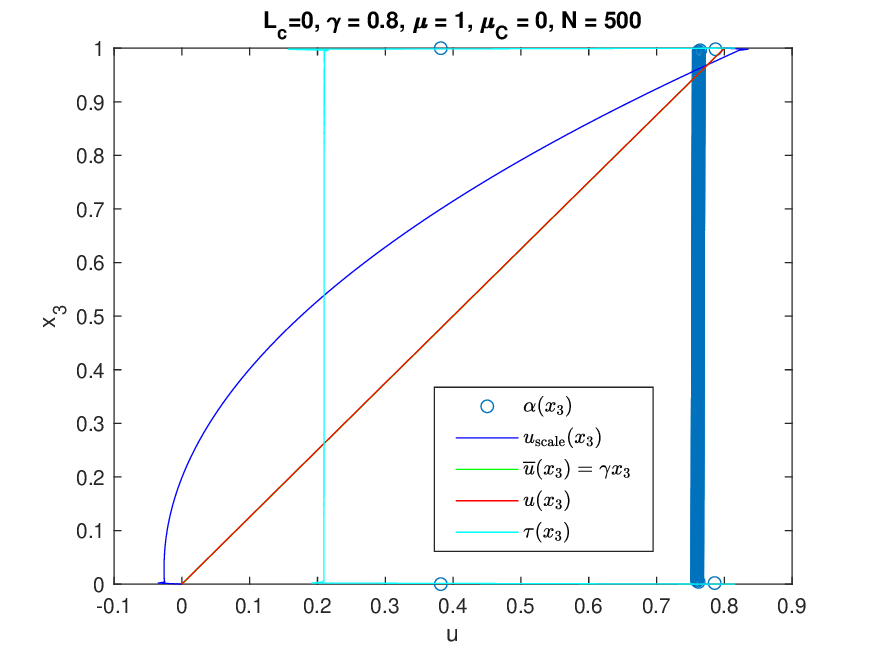,width=6.6cm}}
\put(1.5,0.00){\psfig{figure=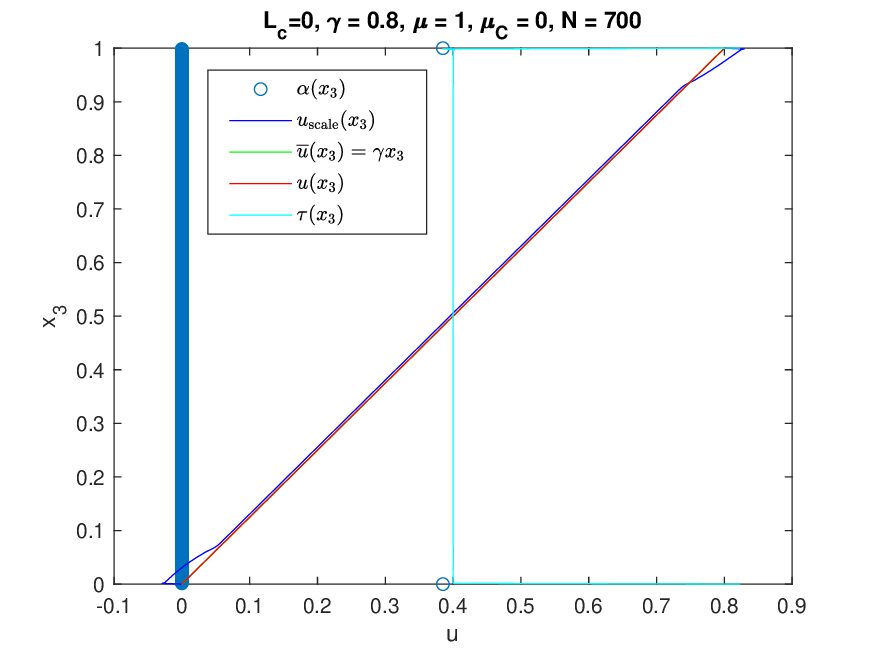,width=6.6cm}}
\put(7.5,0.00){\psfig{figure=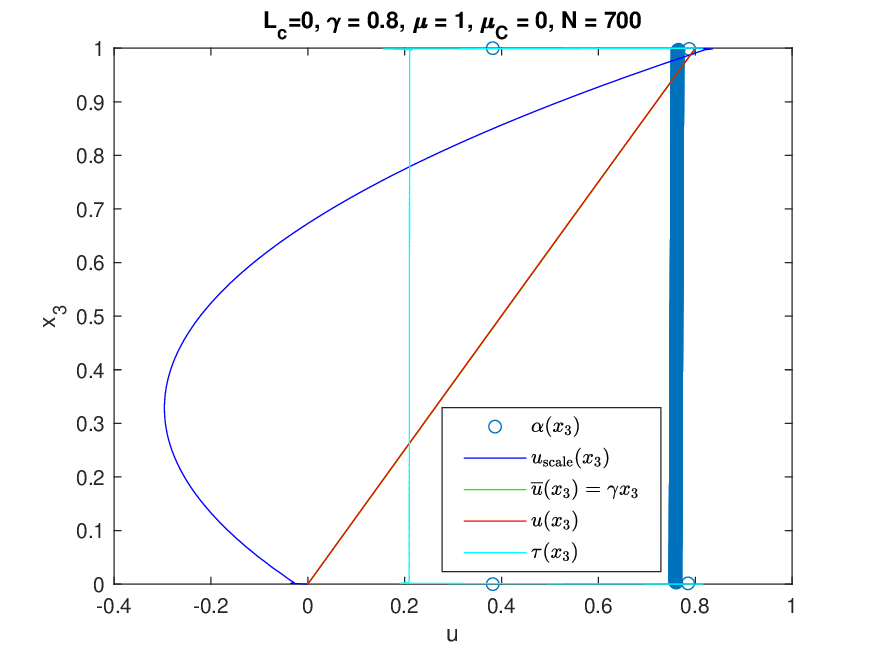,width=6.6cm}}
\end{picture}
\caption{\label{fig8}
Minimizers of $E$ in $\cX^C$ for $\mu=1$, $\mu_c=0$, $\gamma=0.8$ and
$N=49$, $N=149$, $N=500$, $N=700$.
The blue balls are the computed values of $\alpha$ with either
$\alpha\equiv\cblue\alpha_1^-\cn=0$ or
$\alpha\equiv\alpha_3=\eta^{-1}(\gamma)=0.760053$.
In cyan the computed stress tensor $\tau$ which is constant in $\Omega$.
The deformation $u$ is rendered in red, extremely close to
$\overline{u}(x)=\gamma\, x$ rendered in green. The blue line is the rescaled
deformation $u_{\rm scale}:=\gamma\, x+N(u(x)-\gamma\, x)$.}
\end{figure}

\vspace*{2mm}
Fig.~\ref{fig8} displays minimizers of $E$ in $\cX^C$. Except near
$\partial\Omega=\{0,1\}$, the minimizing rotation $\alpha$ takes constant
values with either $\alpha\equiv\cblue\alpha_1^-\cn:=0$ or
$\alpha\equiv\alpha_3:=\eta^{-1}(\gamma)$ in accordance with
Corollary~\ref{cor1}~(ii).
Depending on the initial values at start of the optimization, the BFGS-algorithm
computes one of two different minimizers $u$ which both converge in $\Omega$
to $\overline{u}(x)$ as $N\to\infty$. The first, $u_1$ with
$\alpha\equiv\cblue\alpha_1^-\cn$, is a straight line except near
$\partial\Omega$. The second, $u_2$ with $\alpha\equiv\alpha_3$, leads to a
scaled function $u_\mathrm{scale}(x):=\gamma\, x+N\big(u_2(x)-\gamma\, x\big)$
which is a bended curve, see Fig.~\ref{fig8}. The energy levels of both computed
local minimizers are extremely close with
\begin{equation}
\label{Wval}
E(u_1,\cblue\alpha_1^-\cn)=0.32001744,\qquad\qquad E(u_2,\alpha_3)=0.320018,
\end{equation}
confirming Eqn.~(\ref{Een13}). Yet, numerical optimization favors
$(\uhom,\cblue\alpha_1^-\cn)$ over $(\uhom,\alpha_3)$, underlining that the
energy landscape for geometrically nonlinear Cosserat materials is extremely
complicated and emphasizing why it is so hard to numerically compute
the correct global minimizers.
In comparison to (\ref{Wval}), the theoretical infimal
energy is $E(\overline{u},\cblue\alpha_1^-\cn)=\frac{\mu}{2}\gamma^2=0.32$, but
$(\overline{u},\cblue\alpha_1^-\cn)$ violates the consistent coupling condition
(\ref{cc}).

\begin{figure}[pht]
\unitlength1cm
\begin{picture}(11.5,16.3)
\put(1.5,11.00){\psfig{figure=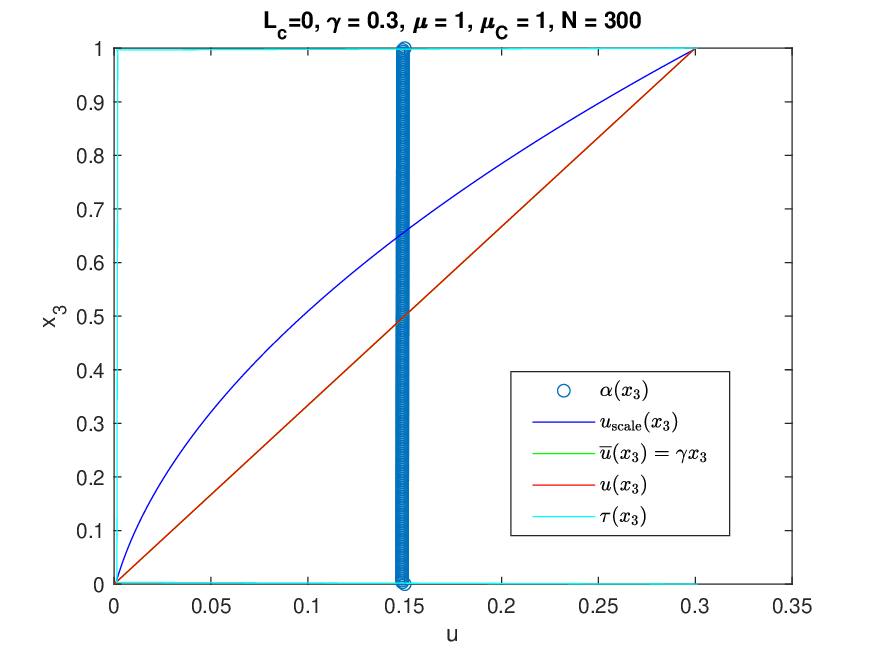,width=6.6cm}}
\put(7.5,11.00){\psfig{figure=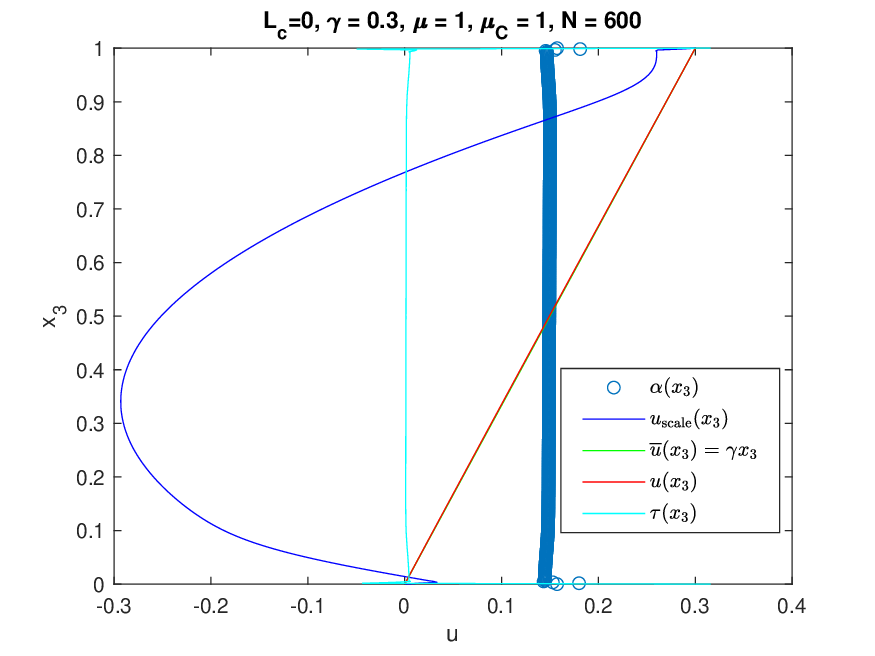,width=6.6cm}}
\put(1.5,5.50){\psfig{figure=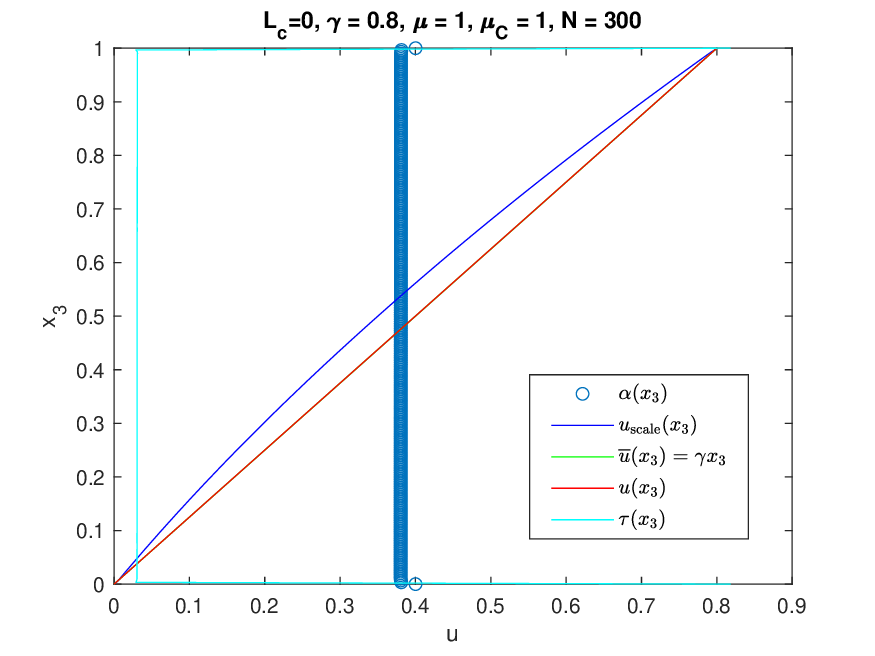,width=6.6cm}}
\put(7.5,5.50){\psfig{figure=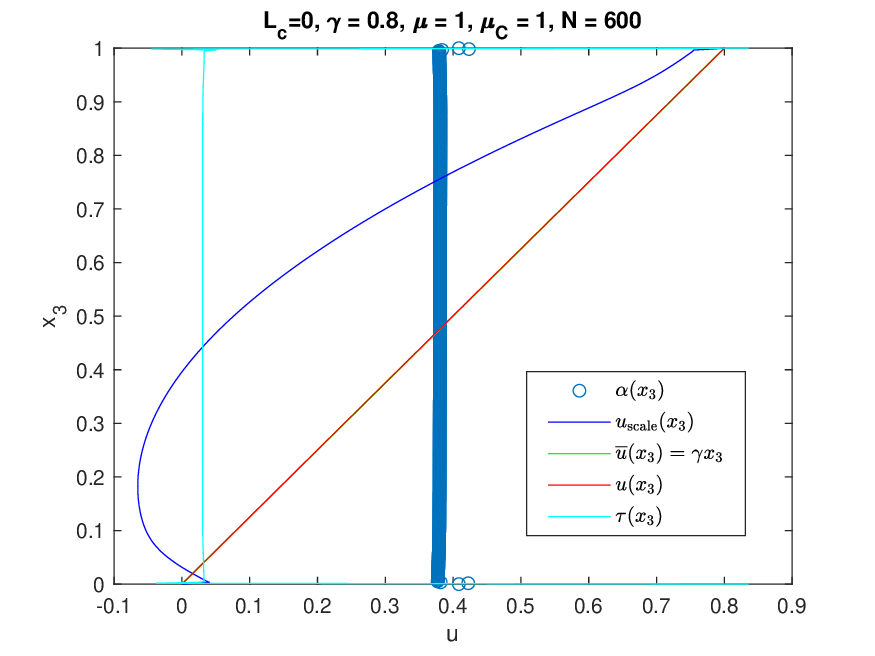,width=6.6cm}}
\put(1.5,0.00){\psfig{figure=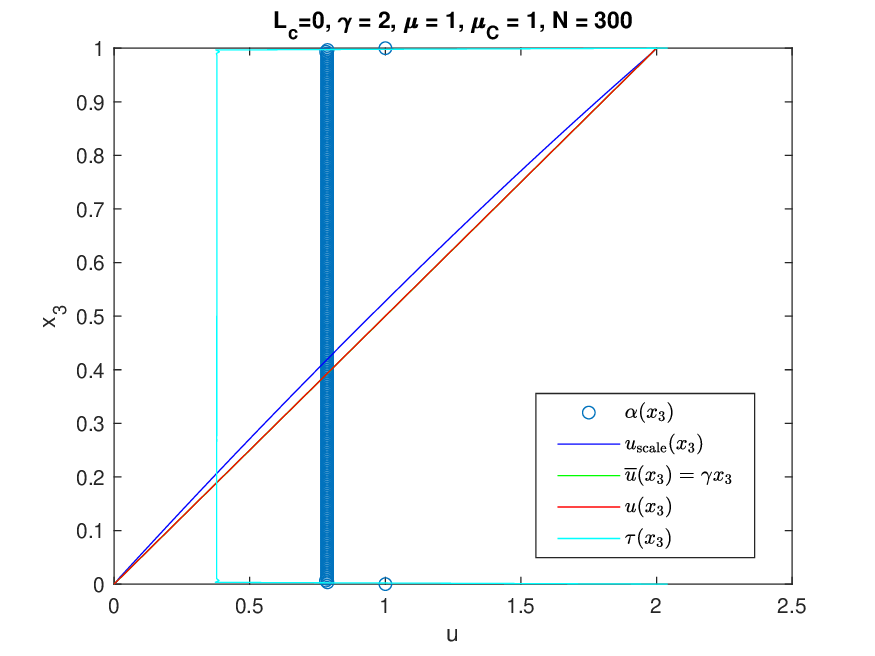,width=6.6cm}}
\put(7.5,0.00){\psfig{figure=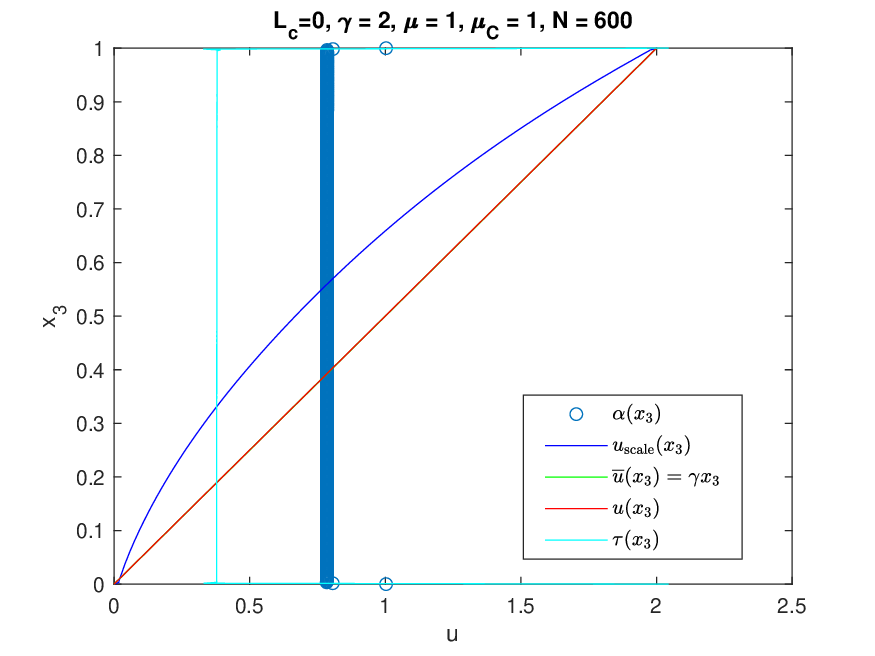,width=6.6cm}}
\end{picture}
\caption{\label{fig9}
Minimizers of $E$ in $\cX^C$ for $\mu=\mu_c=1$ and $N=300$ (left), $N=600$
(right) and $\gamma=0.3$ top line, $\gamma=0.8$ center, and $\gamma=2.0$ bottom
line. The blue balls are the computed values of $\alpha$ with either
$\alpha\equiv0.14868\sim\alpha_2=\arctan(0.15)=0.14889$ (top),
$\alpha\equiv0.3804\sim\alpha_2=\arctan(0.4)=0.380506$ (center) and
$\alpha\equiv0.7852\sim\alpha_2=\arctan(1)=0.785398$ (bottom) as predicted
by Cor.~\ref{cor1} (i). In cyan the computed stress tensor $\tau$.
The deformation $u$ is rendered in red, extremely close to
$\overline{u}(x)=\gamma\, x$ rendered in green. The blue line is the rescaled
deformation $u_{\rm scale}:=\gamma\, x+N(u(x)-\gamma\, x)$.}
\end{figure}

Fig.~\ref{fig9} studies minimizers in $\cX^C$ related to Corollary~\ref{cor1},
(i). Again, $u_N\to\overline{u}$ in $\Omega$ for $N\to\infty$.
The minimizing microrotation angle $\alpha$ is constant (except near
$\partial\Omega$) and takes the values predicted by Cor.~\ref{cor1}.

\begin{figure}[pht]
\unitlength1cm
\begin{picture}(11.5,11.3)
\put(1.5,5.50){\psfig{figure=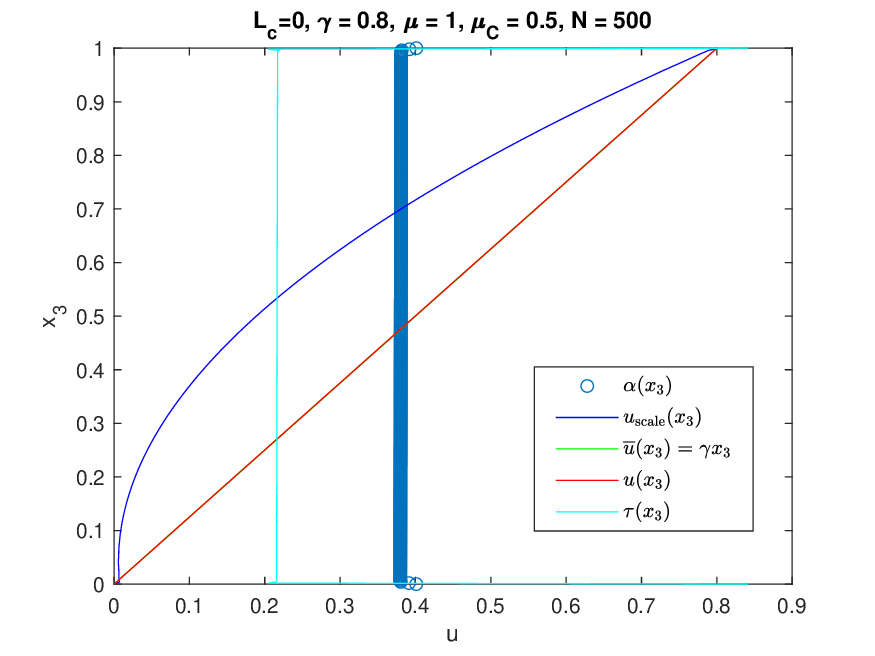,width=6.6cm}}
\put(7.5,5.50){\psfig{figure=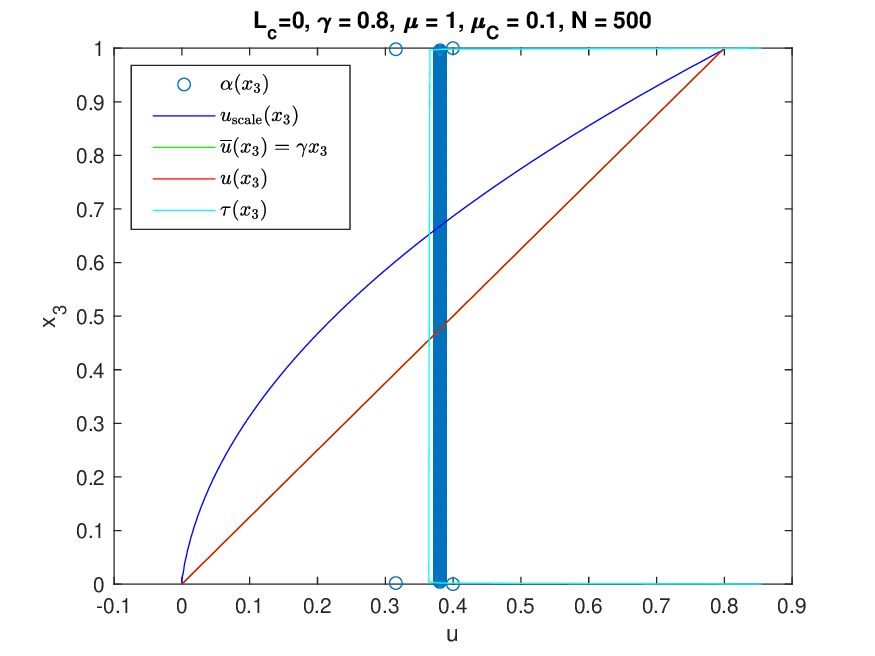,width=6.6cm}}
\put(1.5,0.00){\psfig{figure=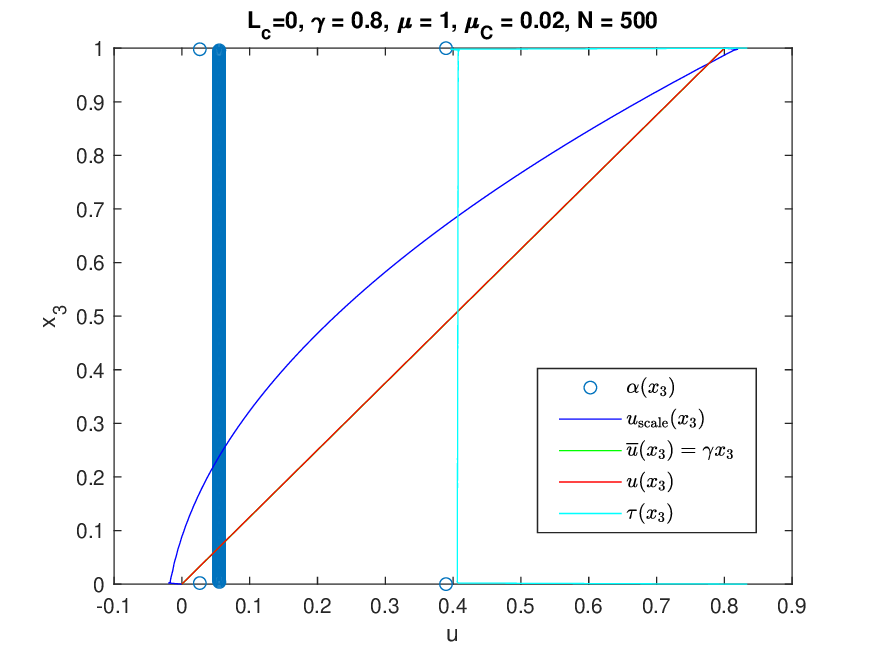,width=6.6cm}}
\put(7.5,0.00){\psfig{figure=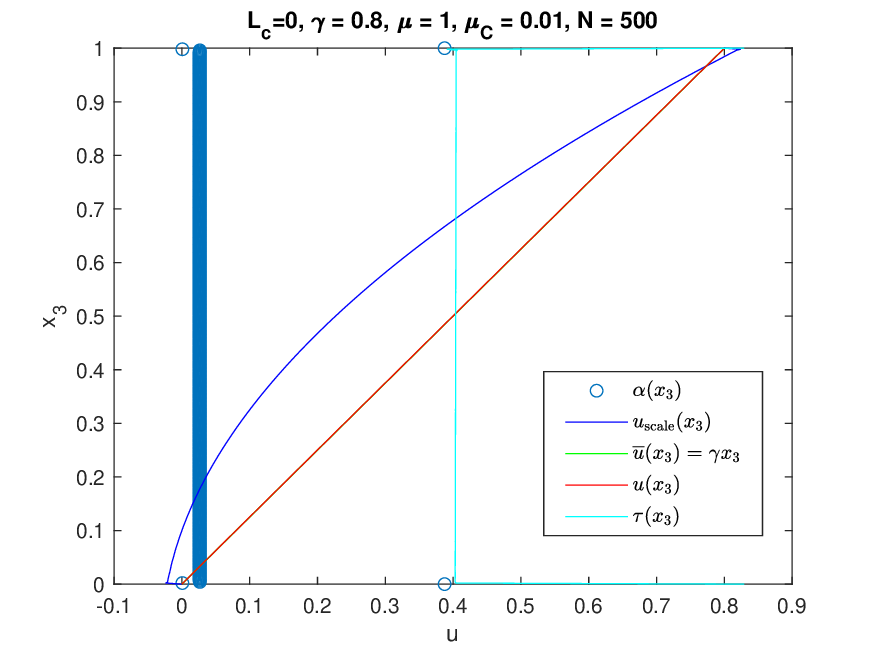,width=6.6cm}}
\end{picture}
\caption{\label{fig10}
Minimizers of $E$ in $\cX^C$ for $\mu=1$, $\gamma=0.8$ and $N=500$.
Top left: $\mu_c=0.5$ with
$\alpha\equiv0.3801\approx\alpha_2=\arctan(\gamma/2)$.
Top right: $\mu_c=0.1$ with $\alpha\equiv0.3804\approx\alpha_2$.
Bottom left: $\mu_c=0.02$ with $\alpha\equiv0.0548\cblue\approx\alpha_1^-\cn$.
Bottom right: $\mu_c=0.01$ with $\alpha\equiv0.026\cblue\approx\alpha_1^-\cn$.
The results underline the validity of Cor.~\ref{cor1} (iii) and demonstrate
that \cblue the value of $\mu_c$ selects either $\alpha_1^\pm$ or $\alpha_2$
as local minimizer of $E(\uhom,\cdot)$, \cn
cf. Remark~\ref{rem7}. The critical value of $\mu_c$ predicted by
(\ref{muccrit}) is $\mu_c^{\rm crit}\approx0.0715$.
In cyan the computed stress tensor $\tau$.
The deformation $u$ is rendered in red, extremely close to
$\overline{u}(x)=\gamma\, x$ rendered in green. The blue line is the rescaled
deformation $u_{\rm scale}:=\gamma\, x+N(u(x)-\gamma\, x)$.}
\end{figure}

Fig.~\ref{fig10} confirms the validity of Corollary~\ref{cor1}, (iii) and
\cblue demonstrate that $(\uhom,\alpha_1^\pm)$ and $(\uhom,\alpha_2)$ are the
local minimizers of $E$ depending on the value of $\mu_c$.\cn



\end{document}